\documentclass[10pt, oneside, reqno]{article}                   				

\makeatletter

	\PassOptionsToPackage{nameinlink,capitalize}{cleveref}
	\PassOptionsToPackage{noend}{algpseudocode}
	\PassOptionsToPackage{breakspaces}{clevethm}
	\usepackage{placeins}
	\usepackage{graphicx}
	\usepackage{xparse}

	\usepackage[%
		arxiv,
		alglinenumber,
		noload={tikz,pgfplots},
		eqreset=section,
		thmreset=section,
	]{myPreamble}

	\let\includetikz\includegraphics
	\graphicspath{{./Pics/Tikz/}}

	\AtBeginDocument{%
		\setcounter{tocdepth}{2}
		\setcounter{secnumdepth}{2}
		\hypersetup{bookmarksdepth = 4}
	}
		\let\oldsubsubsection\subsubsection
		\RenewDocumentCommand{\subsubsection}{s o m}{%
		\IfBooleanTF{#1}
			{%
			\oldsubsubsection*{#3}%
			\addcontentsline{toc}{subsubsection}{\protect\numberline{}#3}%
			}{%
			\IfNoValueTF{#2}
				{\oldsubsubsection{#3}}%
				{\oldsubsubsection[#2]{#3}}%
			}%
		}
	\newlist{claims}{enumerate}{1}
	\setlist[claims,1]{
		label={\it Claim \thedummythm.\oldstylenums{\arabic*}:},
		ref={\thedummythm.\oldstylenums{\arabic*}},
		itemindent=\widthof{\it \thedummythm.Claim \oldstylenums{3}: },
		partopsep=0pt,
		parsep=0pt,
		itemsep=3pt,
		leftmargin=*,
	}
	\newlist{claims*}{enumerate}{1}
	\setlist[claims*,1]{
		label={\it Contradiction claim \oldstylenums{\arabic*}$^*$:},
		ref={\oldstylenums{\arabic*}$^*$},
		partopsep=0pt,
		parsep=0pt,
		itemsep=3pt,
		wide,
		labelindent=0pt,
	}
	\Crefname{claimsi}{Claim}{Claims}

	\Crefname{figure}{Figure}{Figures}
	\newcommand{\secref}{\@ifstar\@@secref\@secref}
	\newcommand{\@secref}[1]{\hyperref[#1]{\S}\ref{#1}}
	\newcommand{\@@secref}[1]{\S\ref*{#1}}

	\crefname{ALG@line}{step}{steps}
	\numberwithin{algorithm}{section}

	\newcommand{\algnamefont}{\sffamily}
	\newcommand{\adaPG}{{\algnamefont adaPGM}}
	\newcommand{\adaPD}{{\algnamefont adaPDM}\@ifnextchar+{\textsuperscript}{}}
	\newcommand{\AdaPG}{{\algnamefont AdaPGM}}
	\newcommand{\AdaPD}{{\algnamefont AdaPDM}\@ifnextchar+{\textsuperscript}{}}

	\newcommand{\refadaPG}{\hyperref[alg:PG]{\adaPG}\@ifstar{}{\ (\cref{alg:PG})}}
	\newcommand{\refAdaPG}{\hyperref[alg:PG]{\AdaPG}\@ifstar{}{\ (\cref{alg:PG})}}

	\newcommand{\refadaPD}{\@ifnextchar+{\refadaPDls}{\@refadaPD}}
	\newcommand{\refadaPDls}[1]{\hyperref[alg:PDls]{\adaPD+}\@ifstar{}{\ (\cref{alg:PDls})}}
	\newcommand{\@refadaPD}{\hyperref[alg:PD]{\adaPD}\@ifstar{}{\ (\cref{alg:PD})}}

	\newcommand{\refAdaPD}{\@ifnextchar+{\refAdaPDls}{\@refAdaPD}}
	\newcommand{\refAdaPDls}[1]{\hyperref[alg:PDls]{\AdaPD+}\@ifstar{}{\ (\cref{alg:PDls})}}
	\newcommand{\@refAdaPD}{\hyperref[alg:PD]{\AdaPD}\@ifstar{}{\ (\cref{alg:PD})}}

	\def\innprod{\@ifstar\@innprod\@@innprod}
	\newcommand{\linop}{\@ifstar{\trans A}{A}}

	\DeclareMathOperator{\diam}{diam}

	\newcommand{\C}{c}
	\renewcommand{\L}{\ell }
	
	\newcommand{\SD}{\mathcal L}
	\newcommand{\U}{\mathcal U}

	\let\OLDdelta\delta     \renewcommand{\delta}{\text{\huge\(\blue\OLDdelta\)}}
	\let\OLDDelta\Delta     \renewcommand{\Delta}{\text{\huge\(\blue\OLDDelta\)}}
	\let\OLDepsilon\epsilon \renewcommand{\epsilon}{\text{\huge\(\blue\OLDepsilon\)}}
	\let\OLDnu\nu           \renewcommand{\nu}{\text{\huge\(\blue\OLDnu\)}}

	\let\indicator\relax
	\DeclareMathOperator{\indicator}{\iota }

	\newcommand{\delksymb}{\OLDdelta }
	\newcommand{\sdcoeff}{\OLDepsilon }
	\newcommand{\V}{\mathcal V}
	\newcommand{\const}{\OLDnu }

		\newcommand{\ck}{\@ifstar{\C_{k+1}}{\C_k}}
		\newcommand{\Lk}{\@ifstar{\L_{k+1}}{\L_k}}
		\newcommand{\Uk}{\@ifstar{\U_{k+1}}{\U_k}}

		\newcommand{\delk}{\@ifstar{\delksymb_{k+1}}{\delksymb_k}}%
		\newcommand{\mytheta}{\vartheta}%
		\newcommand{\gamk}{\@ifstar{\gamma_{k+1}}{\gamma_k}}%
		\newcommand{\Hk}{\@ifstar{H_{k+1}}{H_k}}%
		\newcommand{\rhok}{\@ifstar{\rho_{k+1}}{\rho_k}}%
		\newcommand{\sigk}{\@ifstar{\sigma_{k+1}}{\sigma_k}}%
		\newcommand{\xik}{\@ifstar{\xi_{k+1}}{\xi_k}}%

		\newcommand{\Fwk}{\@ifstar\@@Fwk\@Fwk}
		\newcommand{\Fw@twoargs}[2][]{{\ifstrempty{#2}{\id-\gamma_{#1}\nabla f}{#2-\gamma_{#1}\nabla f(#2)}}}
		\newcommand{\@Fwk}[2][k]{\Fw@twoargs[#1]{#2}}
		\newcommand{\@@Fwk}[2][k+1]{\Fw@twoargs[#1]{#2}}

		\newcommand{\FBk}{\@ifstar\@@FBk\@FBk}
		\newcommand{\FB@twoargs}[2][]{{\prox_{\gamma_{#1}g}(\Fw@twoargs[#1]{#2})}}
		\newcommand{\@FBk}[2][k]{\FB@twoargs[#1]{#2}}
		\newcommand{\@@FBk}[2][k+1]{\FB@twoargs[#1]{#2}}

		\newcommand{\PDx}{\@ifstar\@@PDx\@PDx}
		\newcommand{\PDx@threeargs}[3][]{{\prox_{\gamma_{#1}g}\left(\Fw@twoargs[#1]{#2}-\gamma_{#1}\linop*\ifstrempty{#3}{{}\cdot{}}{#3}\right)}}
		\newcommand{\@PDx}[3][k]{\PDx@threeargs[#1]{#2}{#3}}
		\newcommand{\@@PDx}[3][k+1]{\PDx@threeargs[#1]{#2}{#3}}

\makeatother

\newcommand{\TheTitle}{Adaptive proximal algorithms for convex optimization under local Lipschitz continuity of the gradient}

\newcommand{\TheFunding}{%
	This work was supported by:
	the Research Foundation Flanders (FWO) postdoctoral grant 12Y7622N and research projects G081222N, G033822N, and G0A0920N;
	Research Council KU Leuven C1 project No. C14/18/068;
	European Union's Horizon 2020 research and innovation programme under the Marie Skłodowska-Curie grant agreement No. 953348;
	Japan Society for the Promotion of Science (JSPS) KAKENHI grant JP21K17710.%
}
\newcommand{\TheKeywords}{%
	Convex minimization
	\(\cdot\)
	proximal gradient method
	\(\cdot\)
	primal-dual algorithms
	\(\cdot\)
	locally Lipschitz gradient
	\(\cdot\)
	linesearch-free adaptive stepsizes%
}
\newcommand{\TheSubjclass}{%
	65K05
	\(\cdot\)
	90C06
	\(\cdot\)
	90C25
	\(\cdot\)
	90C30
	\(\cdot\)
	90C47%
}
	\title{\TheTitle\thanks{\TheFunding}}
	\author{%
		Puya Latafat\thanks{%
			\protect\TheAddressKU.
			{\it E-mails:} \sf
			\{%
				\emailLink[puya.latafat@kuleuven.be]{puya.latafat}%
			,%
				\emailLink[panos.patrinos@kuleuven.be]{panos.patrinos}%
			\}%
			\emailLink[puya.latafat@kuleuven.be,panos.patrinos@kuleuven.be]{@kuleuven.be}%
		}%
		\and
		Andreas Themelis\thanks{%
			\protect\TheAddressKUJ.
			{\it E-mail:} \sf\emailLink{andreas.themelis@ees.kyushu-u.ac.jp}%
		}%
		\and
		Lorenzo Stella\thanks{%
			\protect\TheAddressAmazon\emph{ (work done prior to joining Amazon)}.
			{\it E-mail:} \sf\emailLink{lorenzostella@gmail.com}%
		}%
		\and
		Panagiotis Patrinos\footnotemark[2]%
	}

	\date{}

\begin{document}

	\maketitle

	\begin{abstract}
		Backtracking linesearch is the de facto approach for minimizing continuously differentiable functions with locally Lipschitz gradient.
		In recent years, it has been shown that in the convex setting it is possible to avoid linesearch altogether, and to allow the stepsize to adapt based on a local smoothness estimate without any backtracks or evaluations of the function value.
		In this work we propose an adaptive proximal gradient method, \refadaPG*, that uses novel estimates of the local smoothness modulus which leads to less conservative stepsize updates and that can additionally cope with nonsmooth terms.
		This idea is extended to the primal-dual setting where an adaptive three-term primal-dual algorithm, \refadaPD*, is proposed which can be viewed as an extension of the {\algnamefont PDHG} method.
		Moreover, in this setting the ``essentially'' fully adaptive variant \refadaPD+* is proposed that avoids evaluating the linear operator norm by invoking a backtracking procedure, that, remarkably, does not require extra gradient evaluations.
		Numerical simulations demonstrate the effectiveness of the proposed algorithms compared to the state of the art.
	\end{abstract}

	\AMS{\TheSubjclass}

	\tableofcontents

	\section{Introduction}\label{sec:introduction}
		Backtracking linesearch is one of the most successful ideas in smooth optimization.
		It is well known that gradient descent with linesearch converges under mild differentiability assumptions \cite[\S1.2]{bertsekas2016nonlinear}.
		Even under Lipschitz gradient continuity such techniques can lead to significant speedups compared to using a constant stepsize dictated by a global Lipschitz modulus due to their ability to adapt to the local geometry of the problem.
		In this work we explore an alternative approach which can cope with nonsmooth formulations, does not require any backtracking procedures or function value evaluations, and yet only requires \emph{local} Lipschitz gradient continuity of the differentiable term.
		The key property that allows for this improvement is the assumption of convexity.

		Our work is inspired by \cite{malitsky2020adaptive} where an adaptive gradient method \(x^{k+1} = \Fwk*{x^k}\) is studied with stepsizes updated by the rule
		\begin{equation}\label{eq:Malitsky:gamk}
			\gamk*
		=
			\min\set{
				\gamk\sqrt{1+\tfrac{\gamk}{\gamma_{k-1}}},\,
				\frac{1}{2L_k}
			},
		\end{equation}
		where
		\begin{equation}\label{eq:MML_k}
			L_k
		\coloneqq
			\frac{
				\|\nabla f(x^k) - \nabla f(x^{k-1})\|
			}{
				\|x^k - x^{k-1}\|
			}
		\end{equation}
		is a local Lipschitz estimate of \(\nabla f\).
		The idea of using an estimate for the Lipschitz modulus has also been explored in the setting of variational inequalities \cite{yang2018modified,thong2020self,bot2023relaxed,Yang2021Self,bohm2022solving}, but often at the cost of enforcing the stepsize sequence to be nonincreasing, which can lead to slow convergence.
		Allowing the stepsize to increase is a crucial feature of \eqref{eq:Malitsky:gamk} which our proposed methods maintain.
		It is worth noting that in \cite{malitsky2020golden} another adaptive scheme, {\algnamefont aGRAAL}, was proposed for hemivariational inequalities which also allows for increasing stepsizes (see also \cite{alacaoglu2023beyond}).
		Recently, in the setting of the gradient method for Lipschitz smooth minimization, \cite{grimmer2023accelerated} advanced an interesting choice of stepsizes according to predefined cyclic patterns. A similar idea appears in \cite{altschuler2023acceleration} that adopts nonrepeating \emph{fractal-like} patterns.
		These methods provably yield improved worst-case rates over the standard gradient method but are bound to globally Lipschitz smooth problems. Moreover,  they incorporate predetermined stepsize sequences agnostic to the local geometry of the cost function.
		Interestingly, \refadaPG*, the adaptive scheme presented in this paper, as well as the one in \cite{malitsky2020adaptive}, automatically leads to sequences of large stepsizes that exhibit a seemingly cyclic behavior, see \cref{fig:gam} and the discussion in \cref{sec:observations}.
		We also mention recent works \cite{teboulle2022elementary,marumo2022parameter} whose adaptive rules are designed to guarantee worst-case rates, and \cite{attouch2023fast} that exploits a continuous-time viewpoint to develop adaptive algorithms.

		In \cite{malitsky2020adaptive} it was observed that the line of proof therein does not provide any route for generalization to the composite proximal setting.
		Additionally to showing that this is in fact possible, in this work we will actually provide larger stepsizes.
		As better detailed in the discussion before \cref{sec:CL}, the improvement is partially due to tighter estimates of the geometry of \(f\) in which both Lipschitz and (inverse) cocoercivity estimates are taken into account, namely
		\begin{subequations}\label{eq:CL}
			\begin{align}
				\Lk
			\coloneqq{} &
				\frac{
					\innprod{\nabla f(x^{k-1})-\nabla f(x^k)}{x^{k-1}-x^k}
				}{
					\|x^{k-1}-x^k\|^2
				}
			\shortintertext{and}
				\ck
			\coloneqq{} &
				\frac{
					\|\nabla f(x^{k-1})-\nabla f(x^k)\|^2
				}{
					\innprod{\nabla f(x^{k-1})-\nabla f(x^k)}{x^{k-1}-x^k}
				}
			\end{align}
		\end{subequations}
		for a pair of points \(x^{k-1},x^k\in\R^n\).
		Noting that the denominator of \(\Lk\) or \(\ck\) is zero iff \(\nabla f(x^k)-\nabla f(x^{k-1})=0\) (the latter owing to the Baillon-Haddad theorem \cite[Cor. 10]{baillon1977quelques}; see \cref{lem:enlarged:ck} for the local version needed here), we stick to the convention \(\frac00=0\) so that both \(\Lk\) and \(\ck\) are (well-defined, positive) real numbers.
		Throughout, we shall also adhere to \(\frac10=\infty\).
		Note that \(\Lk\) and \(\ck\) are the inverse of the Barzilai-Borwein stepsize choices \cite{barzilai1988two}, which have been considered in the setting of gradient descent \cite{raydan1993barzilai,dai2005projected} but whose convergence results are limited to the quadratic setting (see also \cite{tan2016barzilai,burdakov2019stabilized} for extensions).
		Our proposed adaptive proximal gradient scheme \refadaPG{} combines the two estimates and involves the update rule
		\begin{equation}\label{eq:PG:gamk}
			\gamk*
		=
			\min\set{
				\gamk\sqrt{1+\tfrac{\gamk}{\gamma_{k-1}}},\,
				\frac{
					\gamk
				}{
					2\sqrt{\left[\gamk\Lk(\gamk\ck-1)\right]_+}
				}
			}
		\end{equation}
		on the stepsize.
		Note that whenever \(\gamk\ck \leq 1\) the update reduces to \(\gamk*=\gamk\sqrt{1+\tfrac{\gamk}{\gamma_{k-1}}}\), effectively strictly increasing the stepsize.
		Regardless, this update is easily seen to be less conservative than \eqref{eq:Malitsky:gamk}, the one prescribed in \cite[Alg. 1]{malitsky2020adaptive}; see \cref{rem:PG:Malitsky} for the details.
		We also point out the recent follow-up work \cite{malitsky2023adaptive} of \cite{malitsky2020adaptive}, subsequent to the preprint version of our manuscript, that also considers the proximal gradient setting (with a different stepsize update) and improves the second term in \eqref{eq:Malitsky:gamk} by a factor of \(\sqrt 2\) in the smooth case.

		In the second part of the paper this idea is extended to the primal-dual setting to address more general problems of the form
		\begin{equation}\label{eq:PD}
			\minimize_{x\in\R^n}\; \varphi(x) \coloneqq f(x) + g(x) + h(\linop x),
		\end{equation}
		where \(\linop\) is a linear mapping, \(g\) and \(h\) are (possibly nonsmooth) extended-real-valued convex functions, and \(f\) is a convex function typically assumed to have Lipschitz continuous gradient (this is relaxed to local Lipschitz continuity here, cf. \cref{ass:PD}).

		In the past decade primal-dual algorithms have gained a lot of popularity in areas ranging from machine learning and signal processing to control \cite{combettes2011proximal,Sra2012optimization,komodakis2015playing,jezierska2012primal,latafat2018plug,latafat2019new}.
		Their popularity is primarily due to their ability to achieve \emph{full splitting} on composite problems of the form \eqref{eq:PD}.
		Moreover, inherent properties of first-order operations facilitate block-coordinate and distributed variants, see for instance \cite{bianchi2014primal,latafat2016new,fercoq2019coordinate,latafat2022block,latafat2022bregman} and the references therein.

		There is a large body of literature on primal-dual algorithms; see, e.g., \cite{chambolle2011first,drori2015simple,condat2013primal,bot2013douglas,latafat2019new,yan2018new}.
		Despite employing different techniques in their convergence analysis, the majority of existing methods rely on establishing a Fej\'er-type inequality.
		In fact, most can be viewed as intelligent applications of a monotone splitting technique such as forward-backward, Douglas-Rachford, and forward-backward-forward splittings for solving the associated primal-dual optimality conditions; see, \eg \cite{he2012convergence,vu2013splitting,condat2013primal,bot2013douglas,combettes2012primal}.
		More recently, the introduction of new splitting techniques such as
		{\algnamefont AFBA} \cite{latafat2017asymmetric,latafat2018primal},
		{\algnamefont NOFOB} \cite{giselsson2021nonlinear},
		forward-Douglas-Rachford-forward \cite{ryu2020finding},
		forward-backward-half forward  \cite{bricenoarias2018forward},
		forward-reflected-backward \cite{malitsky2020forward},
		has led to new primal-dual algorithms.
		We remark also that when \(\linop=\id\) in \eqref{eq:PD} one can directly solve the problem without any lifting by using the three-term splitting \cite{davis2017three}.
		There exists also an adaptive/linesearch variant of this algorithm (see, e.g., \cite{pedregosa2018adaptive}) which however requires potentially costly extra gradient evaluations during the backtracking procedure.

		Unlike the above described correspondence with splitting techniques, our proposed method cannot be viewed as an instance of any splitting technique for solving \emph{general} monotone inclusions, in that it relies heavily on the knowledge that the operators involved are subdifferentials of convex functions.
		Although the proposed idea is extendable to other primal-dual methods such as those in \cite{latafat2017asymmetric}, our focus here is on an adaptive variant of {\algnamefont PDHG} \cite{chambolle2011first,vu2013splitting,condat2013primal}.
		An interesting algorithm in this line of work was proposed in \cite{vladarean2021first} which however cannot handle the third nonsmooth term \(g\) in \eqref{eq:PD}.
		In \refadaPD\ we provide a different stepsize rule that not only can handle \eqref{eq:PD} but also inherits the same idea of using tighter estimates \(\ck, \Lk\) as in the case of the proximal gradient method.

		A second consideration for primal-dual methods is that in the usual (nonadaptive) setting the primal-dual stepsizes \(\gamma, \sigma\) should typically satisfy a condition of the form \(\gamma \sigma \|\linop\|^2 \leq 1 - \tfrac{\gamma L_f}2\) (see for instance \cite[Thm. 3.1]{condat2013primal} and \cite[Prop. 5.1]{latafat2017asymmetric}).
		In practical applications the norm of the linear operator may be costly to compute and can lead to smaller stepsizes, and thus slower convergence.
		Recently, in the setting where \(f\equiv 0\), a linesearch procedure was proposed in \cite{malitsky2018first} for selecting the stepsizes based on an estimate of the norm of the linear operator.
		A linesearch extension of \cite{malitsky2020golden} in the primal-dual setting was also proposed in \cite{chang2022golden}.
		Instead, the linesearch procedure that we propose naturally integrates our adaptive primal-dual algorithm \refadaPD+ to handle the more general problem \eqref{eq:PD} without any extra gradient evaluations during the backtracks.

		\subsection{Contributions}%
			The main contributions of the paper are summarized below.
			\begin{enumerate}[%
				leftmargin=0pt,
				label={\arabic*.},
				align=left,
				labelwidth=1em,
				itemindent=\labelwidth+\labelsep,
			]
			\item
				We propose a nonmonotone adaptive stepsize rule for the proximal gradient method that departs from the usual linesearch technique.
				In contrast to backtracking linesearch, the new approach eliminates the need for backtracks or function value evaluations altogether.
				More importantly, the proposed algorithm does not require any parameter tuning and can quickly recover from a bad stepsize initialization.
				This is achieved by adapting the stepsize to the local geometry of the smooth function, combining local estimates of cocoercivity and Lipschitz moduli of the differentiable term along the last two iterates.
				Compared to \cite{malitsky2020adaptive}, even when restricted to the case of gradient descent, the proposed approach allows for less restrictive stepsizes.
				Through this observation, convergence of the aforementioned work in the proximal case follows immediately as a by-product of our analysis.
			\item
				This idea is extended to the primal-dual setting where an adaptive three-term splitting for composite minimization problems is developed.
				The proposed algorithm can be viewed as an extension of the Condat-V\~u algorithm \cite{condat2013primal,vu2013splitting}, which in turn is an extension of the {\algnamefont PDHG} algorithm \cite{chambolle2011first}.
			\item
				As a final contribution, an ``essentially'' fully adaptive variant of the primal-dual method is presented.
				This is meant in the sense that it no longer requires evaluating the norm of the linear operator \(\linop\), and is thus ``fully'' adaptive, but only ``essentially'' so, for all this comes at the expense of performing a backtracking to potentially correct the given (local) estimates.
				Remarkably nevertheless, the proposed linesearch does not require any extra gradient evaluations and can thus be implemented efficiently.
			\end{enumerate}

		\subsection{Organization}%
			We conclude this section by introducing the adopted notation.
			The proposed adaptive proximal gradient method \refadaPG* is formally studied in \cref{sec:PG}.
			The underlying idea is then extended to the primal-dual setting in \cref{sec:PD}, where \refadaPD* is presented that can handle one additional nonsmooth term composed with a linear operator.
			The issue of estimating the norm of the linear operator is resolved through the introduction of a linesearch procedure in \cref{sec:PDls}.
			The convergence results for both variants of the primal-dual algorithm are presented in a unified fashion in \cref{sec:convergence:PD} with some of the proofs deferred to
			\cref{sec:PD:convergence}.
			Numerical simulations for the proposed algorithms are presented in \cref{sec:simulations}, together with a commentary on some empirical observations.
			\Cref{sec:conclusions} concludes the paper.

		\subsection{Notation}\label{sec:notation}%
			The set of real and extended-real numbers are \(\R\coloneqq(-\infty,\infty)\) and \(\Rinf\coloneqq\R\cup\set\infty\), while the positive and strictly positive reals are \(\R_+\coloneqq[0,\infty)\) and \(\R_{++}\coloneqq(0,\infty)\).
			We use the notation \([x]_+ = \max\set{0, x}\).
			With \(\id\) we indicate the identity function defined on a suitable space.
			We denote by \(\innprod{{}\cdot{}}{{}\cdot{}}\) and \(\|{}\cdot{}\|\) the standard Euclidean inner product and the induced norm.
			Given a set \(\mathcal D\subseteq\R^n\), with \(\interior \mathcal D\), \(\relint \mathcal D\) and \(\boundary \mathcal D\) we respectively denote its interior, relative interior, and boundary, and for a sequence \(\seq{x^k}\) we write \(\seq{x^k}\subseteq \mathcal D\) to indicate that \(x^k\in \mathcal D\) for all \(k\in\N\).
			The diameter of \(\mathcal D\) is \(\diam\mathcal D\coloneqq\sup\set{\|x-y\|}[x,y\in\mathcal D]\), while its indicator function is \(\indicator_{\mathcal D}\), namely \(\indicator_{\mathcal D}(x)=0\) if \(x\in \mathcal D\) and \(\infty\) otherwise.
			Given two sets \(\mathcal D,\mathcal D'\subseteq\R^n\), the notation \(\mathcal D+\mathcal D'\coloneqq\set{x+x'}[x\in\mathcal D,\ x'\in\mathcal D']\) indicates their Minkowski sum.

			The notation \(\cball{\bar x}{r}\coloneqq\set{x}[\|x-\bar x\|\leq r]\) indicates the closed ball centered at \(\bar x\) and with radius \(r\).
			The \DEF{domain} and \DEF{epigraph} of an extended-real-valued function \(\func{h}{\R^n}{\Rinf}\) are the sets
			\(
				\dom h
			\coloneqq
				\set{x\in\R^n}[
					h(x)<\infty
				]
			\)
			and
			\(
				\epi h
			\coloneqq
				\set{(x, c )\in\R^n\times\R}[
					h(x)\leq c
				]
			\).
			Function \(h\) is said to be \DEF{proper} if \(\dom h\neq\emptyset\), and \DEF{lower semicontinuous (lsc)} if \(\epi h\) is a closed subset of \(\R^{n+1}\).
			We say that \(h\) is \emph{level bounded} if its \( c \)-sublevel set
			\(
				\lev_{\leq c }h
			\coloneqq
				\set{x\in\R^n}[
					h(x)\leq c
				]
			\)
			is bounded for all \(c\in\R\).
			The \DEF{conjugate} of \(h\) is defined by
			\(
				\conj h(y)
			\coloneqq
				\sup_{x\in\R^n}\set{\innprod yx-h(x)}
			\).

			We say that a differentiable function \(\func{h}{\R^n}{\R}\) has \DEF{locally Lipschitz continuous gradient} if for every convex and compact set \(\V\subset\R^n\) there exists \(L_{h,\V}>0\) such that
			\(
				\|\nabla h(x)-\nabla h(x')\|\leq L_{h,\V}\|x-x'\|
			\)
			holds for all \(x,x'\in\mathcal\V\).\footnote{%
				By virtue of \cite[Thm. 9.2]{rockafellar2009variational}, this condition is equivalent to \DEF{strict continuity} of \(\nabla h\) at every point, in the sense of \cite[Def. 9.1]{rockafellar2009variational}.
			}

	\section{Adaptive proximal gradient method}\label{sec:PG}
		The proximal gradient method ({\algnamefont PGM}) is the natural extension of gradient descent for constrained and nonsmooth problems.
		It addresses nonsmooth minimization problems by splitting them into the sum of two terms as follows:
		\begin{equation}\label{eq:PG}
			\minimize_{x\in\R^n}\; \varphi(x) \coloneqq f(x)+g(x).
		\end{equation}
		Throughout this section the following underlying assumptions are imposed on problem \eqref{eq:PG}.

		\begin{assumption}[Requirements for the proximal gradient setting]\label{ass:PG}
			The following hold in problem \eqref{eq:PG}:
			\begin{enumeratass}
			\item \label{ass:f}%
				\(\func{f}{\R^n}{\R}\) is convex and has locally Lipschitz continuous gradient.
			\item
				\(\func{g}{\R^n}{\Rinf}\) is proper, lsc, and convex.
			\item
				A solution exists: \(\argmin\varphi\neq\emptyset\).
			\end{enumeratass}
		\end{assumption}

		In addition to the gradient of the differentiable term, the fundamental oracle of {\algnamefont PGM} is the \emph{proximal mapping} \cite[Def. 6.1]{beck2017first}
		defined as
		\[
			\prox_{\tau g}(x)
		\coloneqq
			\argmin_{w\in\R^n}
			\set{g(w) + \tfrac{1}{2\tau}\|w-x\|^2},
		\]
		where \(\tau>0\) is a given stepsize.
		In the convex setting the proximal map is single valued and in fact Lipschitz continuous.
		It is well known that for many applications of interest such as constrained or regularized problems the nonsmooth term admits closed form proximal operator (e.g., projection on sets, shrinkage operator, etc.).
		The most common variant of {\algnamefont PGM} involves a constant stepsize that is upper bounded by \(\nicefrac{2}{L_f}\), where \(L_f\) is the global Lipschitz constant of \(\nabla f\).
		A common strategy in practice is to estimate such modulus via backtracking linesearch.

		\begin{algorithm}[t]
			\caption{Adaptive proximal gradient method (\protect\adaPG)}
			\label{alg:PG}
			\begin{algorithmic}[1]
			\itemsep=3pt
			\Require
				starting point \(x^{-1}\in\R^n\)
				~and~
				stepsizes \(\gamma_0\geq\gamma_{-1}>0\)
			\Initialize
				\(x^0=\FBk[0]{x^{-1}}\)
			\item[\algfont{Repeat for} \(k=0,1,\ldots\) until convergence]
			\State \label{state:PG:gamk}%
				With \(\Lk\) and \(\ck\) as in \eqref{eq:CL}, define the stepsize as
			\hfill{\gray\footnotesize\sf // \([r]_+\coloneqq\max\set{r,0}\)}
				\[
					\gamk*
				=
					\gamk\min\set{
						\sqrt{1+\tfrac{\gamk}{\gamma_{k-1}}},\,
						\frac{
							1
						}{
							2\sqrt{\left[\gamk\Lk(\gamk\ck-1)\right]_+}
						}
					}
				\]
			\State
				\(x^{k+1}=\FBk*{x^k}\)
			\end{algorithmic}
		\end{algorithm}

		\subsection{Algorithmic overview}\label{sec:PG:overview}
			Diverging from the linesearch technique, we propose \refadaPG* that adaptively selects the stepsizes based on the estimates \(\ck, \Lk\) as in \eqref{eq:CL}.
			Note that an adaptive gradient method with \(g\equiv0\) was proposed in \cite{malitsky2020adaptive} with a different update rule, see \eqref{eq:Malitsky:gamk}, where it was observed that the adopted line of proof \emph{does not seem to provide any route for generalization} to account for a nonsmooth term.
			This appears to be fundamentally due to the fact that the analysis therein revolves around the fact that the difference of consecutive iterates is a multiple of the gradient.
			In contrast, we circumvent this by combining the subgradient inequality for the nonsmooth term \(g\) at \emph{three} different pairs of points, namely \((x^\star, x^{k+1})\), \((x^{k+1}, x^k)\), and \((x^{k-1}, x^k)\).
			In addition, the combined use of the quantities \(\Lk\) and \(\ck\) as in \eqref{eq:CL} allows for estimating, along with that of the gradient, the local Lipschitz constant of the forward operator \(\id-\gamk\nabla f\).
			This appears to be fundamental for recovering the update rule \eqref{eq:Malitsky:gamk} of \cite{malitsky2020adaptive}, and in fact leads to the less conservative update rule of \refadaPG*.

			\subsubsection{Initialization and practical considerations}\label{sec:PG:init}%
				\refAdaPG* has the ability to recover from a small stepsizes which may be due to bad initialization, or stumbling upon steep/ill-conditioned regions, by linearly increasing the stepsize (by a factor of at least \(\sqrt2\)) until a value proportional to the inverse of a local Lipschitz estimate is attained; see the proof of \cref{thm:PG:gamk} for the details.
				On the other extreme, too large an initial stepsize is corrected in one iteration thanks to the second term in \eqref{eq:PG:gamk}.
				Nevertheless, this event can also result in the very first algorithmic step being pushed far away without control.

				To mitigate such scenarios at initialization, \(\gamma_0\) can be refined by running offline proximal gradient updates;
				starting from the initial point \(x_{-1}\), \(\gamma_0\) can be updated by the inverse of either one of the quantities in \eqref{eq:CL} or \eqref{eq:MML_k} evaluated between \(x^{-1}\) and the prox-grad point.
				If the updated stepsize is orders of magnitude smaller than the original one, the same procedure may be repeated an additional time.
				This procedure is helpful even for problems with globally Lipschitz gradient continuity.
				Once a reasonable \(\gamma_0\) is obtained, we suggest selecting \(\gamma_{-1}\) small enough such that
				\(
					\gamma_0 \sqrt{1+\nicefrac{\gamma_0}{\gamma_{-1}}}
				\geq
					\frac{1}{2{L_0}},
				\)
				ensuring that \(\gamma_1 \geq \frac{1}{2{L_0}}\). It is important to note that this choice of \(\gamma_{-1}\) doesn't affect the convergence results of \cref{thm:PG:convergence}. In fact, with this initialization \cref{thm:PG:gamk} holds with \(k_0 =1\).

		\subsection{Preliminary lemmas}\label{sec:CL}
			Throughout, we will make use of the following shorthand for the forward operator with stepsize \(\gamk\):
			\[
				\Hk\coloneqq\Fwk{}.
			\]
			The subgradient characterization of the proximal step implies that (see \eqref{eq:PG:subgrad})
			\begin{equation}\label{eq:subgrad:varphi}
				\tfrac{1}{\gamk}\bigl(\Hk(x^{k-1})- \Hk(x^k)\bigr)
			\in
				\partial\varphi(x^k).
			\end{equation}
			This quantity plays an important role in our analysis.
			As we are about to show, the combined adoption of the estimates \(\Lk\) and \(\ck\) provides an estimate of the Lipschitz modulus of not only \(\nabla f\) but also the forward operator \(\Hk\).

			\begin{lemma}\label{thm:PG:cLl}%
				Suppose that \cref{ass:PG} holds, and let \(x^{k-1},x^k\in\R^n\).
				Then, with \(L_k\), \(\Lk\), and \(\ck\) as in \eqref{eq:MML_k} and \eqref{eq:CL} the following hold:
				\begin{enumerate}
				\item \label{thm:PG:CL}%
					\(
						\|\nabla f(x^{k-1})-\nabla f(x^k)\|^2
					=
						\ck\Lk\|x^{k-1}-x^k\|^2
					\), that is, \(\ck\Lk=L_k^2\).
				\item \label{thm:PG:H}%
					\(
						\|\Hk(x^{k-1})-\Hk(x^k)\|^2
					=
						\bigl(1-\gamk\Lk(2-\gamk\ck)\bigr)
						\|x^{k-1}-x^k\|^2
					\).
				\item \label{thm:Lk<=Ck}%
					\(
						\Lk
					\leq
						L_k
					\leq
						L_{f,\V}
					\)
					and \(L_k\leq\ck\), where \(L_{f,\V}\) is a Lipschitz modulus for \(\nabla f\) on a compact convex set \(\V\) containing \(x^{k-1}\) and \(x^k\).
				\end{enumerate}
				\begin{proof}
					The first assertion is of trivial verification, and similarly the third one follows from the Cauchy-Schwarz inequality and the Baillon-Haddad theorem \cite[Cor. 10]{baillon1977quelques}, see also \cite[Cor. 18.17]{bauschke2017convex}.
					To conclude, observe that
					\begin{align*}
						\|\Hk(x^{k-1})-\Hk(x^k)\|^2
					={} &
						\|x^{k-1}-x^k\|^2
						+
						\gamk^2
						\|\nabla f(x^{k-1})-\nabla f(x^k)\|^2
					\\
					&
						-
						2\gamk\innprod{\nabla f(x^{k-1})-\nabla f(x^k)}{x^{k-1}-x^k},
					\end{align*}
					from which assertion \itemref{thm:PG:H} follows.
				\end{proof}
			\end{lemma}

			Our convergence analysis will rely on first establishing boundedness of the generated sequence, thereby entailing the existence of a Lipschitz constant \(L_{f,\V}>0\) for \(\nabla f\) on a bounded convex set \(\V\) that contains the iterates.
			It will then follow from \cref{thm:Lk<=Ck} that both \(L_k\) and \(\Lk\) are upper bounded by this quantity, which in turn will be used to show that this modulus provides a lower bound for the stepsize separating it from zero.
			Before that, we show how the combined use of \(\Lk\) and \(\ck\) can be employed to estimate the progress of the iterates generated by {\algnamefont PGM} with arbitrary stepsizes, not necessarily dictated by the update rule of \refadaPG*.
			To simplify the presentation, we introduce the following notation:
			\begin{equation}\label{eq:rhoP}
				\rhok\coloneqq\tfrac{\gamk}{\gamma_{k-1}}
			\quad\text{and}\quad
				P_k\coloneqq\varphi(x^k)-\min\varphi.
			\end{equation}
			\begin{lemma}\label{thm:PG:descent}%
				Suppose \cref{ass:PG} holds, and consider a sequence \(\seq{x^k}\) generated by {\algnamefont PGM} iterations \(x^{k+1}=\FBk*{x^k}\).
				Then%
				\begin{align*}
				&
					\tfrac{1}{2}\|x^{k+1}-x^\star\|^2
					+
					\gamk*(1+\rhok*)P_k
					+
					\tfrac{1}{4}\|x^k-x^{k+1}\|^2
				\\
				\leq{} &
					\tfrac{1}{2}\|x^k-x^\star\|^2
					+
					\rhok*\gamk*P_{k-1}
					-
					\rhok*^2\gamk\Lk \left(
			 			1-\gamk \ck  		\right)
					\|x^{k-1}-x^k\|^2
				\end{align*}
				holds for any \(k\geq1\) and \(x^\star\in\argmin\varphi\), where \(\Lk\) and \(\ck\) are as in \eqref{eq:CL}, and \(\rhok\) and \(P_k\) as in \eqref{eq:rhoP}.%
				\begin{proof}
					The subgradient characterization
					\begin{equation}\label{eq:PG:subgrad}
						\tfrac{\Hk*(x^k)-x^{k+1}}{\gamk*}
					=
						\tfrac{x^k-x^{k+1}}{\gamk*}-\nabla f(x^k)
					\in
						\partial g(x^{k+1})
					\end{equation}
					of \(x^{k+1}=\FBk*{x^k}\) implies that
					{\mathtight[0.95]%
						\begin{align*}
							0
						\leq{} &
							g(x^\star)-g(x^{k+1})
							+
							\innprod{\nabla f(x^k)}{x^\star- x^{k+1}}
							-
							\tfrac{1}{\gamk*}\innprod{x^k-x^{k+1}}{x^\star-x^{k+1}}
						\\
						={} &
							g(x^\star)-g(x^{k+1})
							+
							\underbracket[0.5pt]{
								\innprod{\nabla f(x^k)}{x^\star- x^{k+1}}
							}_{\text{(A)}}
							+
							\tfrac{1}{2\gamk*}\|x^k-x^\star\|^2
							-
							\tfrac{1}{2\gamk*}\|x^{k+1}-x^\star\|^2
							-
							\tfrac{1}{2\gamk*}\|x^k-x^{k+1}\|^2
						\end{align*}
					}%
					holds for any solution \(x^\star\).
					We next proceed to upper bound the term (A) as
					\begin{align*}
						\text{(A)}
					={} &
						\innprod{\nabla f(x^k)}{x^\star-x^k}
						+
						\innprod{\nabla f(x^k)}{x^k-x^{k+1}}
					\\
					={} &
						\innprod{\nabla f(x^k)}{x^\star-x^k}
						+
						\tfrac{1}{\gamk}\innprod{\Hk(x^{k-1})-x^k}{x^{k+1}-x^k}
						+
						\tfrac{1}{\gamk}\innprod{\Hk(x^{k-1})-\Hk(x^k)}{x^k-x^{k+1}}
					\\
					\overrel[\leq]{\eqref{eq:PG:subgrad}}{} &
						\fillwidthof[c]{\innprod{\nabla f(x^k)}{x^\star-x^k}}{
							f(x^\star) - f(x^k)
						}
						+
						\fillwidthof[c]{\tfrac{1}{\gamk}\innprod{\Hk(x^{k-1})-x^k}{x^{k+1}-x^k}}{
							g(x^{k+1})-g(x^k)
						}
						+
						\underbracket[0.5pt]{
							\tfrac{1}{\gamk}\innprod{\Hk(x^{k-1})-\Hk(x^k)}{x^k-x^{k+1}}
						}_{\text{(B)}}.
					\end{align*}
					We bound the term (B) by Young's inequality with parameter \(\varepsilon_{k+1}\) as
					\begin{align*}
						\text{(B)}
					\leq{} &
						\tfrac{\varepsilon_{k+1}}{2\gamk}\|x^k-x^{k+1}\|^2
						+
						\tfrac{1}{2\varepsilon_{k+1}\gamk}\|\Hk(x^{k-1})-\Hk(x^k)\|^2
					\\
					\overrel{\ref{thm:PG:H}}{} &
						\tfrac{\varepsilon_{k+1}}{2\gamk}\|x^k-x^{k+1}\|^2
						+
						\tfrac{1-\gamk\Lk(2-\gamk\ck)}{2\varepsilon_{k+1}\gamk}\|x^{k-1}-x^k\|^2.
						\numberthis\label{eq:B:CL}
					\end{align*}
					Let \(\varphi_\star \coloneqq \min \varphi\).
					The three inequalities combined give
					\begin{align*}
						0
					\leq{} &
						\varphi_\star-\varphi(x^k)
						+
						\tfrac{1}{2\gamk*}\|x^k-x^\star\|^2
						-
						\tfrac{1}{2\gamk*}\|x^{k+1}-x^\star\|^2
					\\
					&
						+
						\left(\tfrac{\varepsilon_{k+1}}{2\gamk}-\tfrac{1}{2\gamk*}\right)\|x^k-x^{k+1}\|^2
						+
						\tfrac{1-\gamk\Lk(2-\gamk\ck)}{2\varepsilon_{k+1}\gamk}\|x^{k-1}-x^k\|^2.
					\end{align*}
					Using again the subgradient \eqref{eq:PG:subgrad} (since \(\partial\varphi=\nabla f+\partial g\)) one has
					\begin{equation}\label{eq:PG:vk}
						v^k
					\coloneqq
						\tfrac{x^{k-1}-x^k}{\gamk}-(\nabla f(x^{k-1})-\nabla f(x^k))
					\in
						\partial\varphi(x^k),
					\end{equation}
					hence, for any \(\mytheta_{k+1}\geq0\),
					\begin{align*}
						0
					\leq{} &
						\mytheta_{k+1}\left(
							\varphi(x^{k-1})-\varphi(x^k)
							-
							\innprod{v^k}{x^{k-1}-x^k}
						\right)
					\\
					={} &
						\mytheta_{k+1}\left(
							\varphi(x^{k-1})-\varphi(x^k)
							-
							\tfrac{1}{\gamk}\|x^k-x^{k-1}\|^2
							+
							\innprod{\nabla f(x^{k-1})-\nabla f(x^k)}{x^{k-1}-x^k}
						\right)
					\\
					={} &
						\mytheta_{k+1}\left(
							\varphi(x^{k-1})-\varphi(x^k)
							-
							\tfrac{1-\gamk\Lk}{\gamk}\|x^k-x^{k-1}\|^2
						\right).
						\numberthis\label{eq:PG:subgradineq}
					\end{align*}
					By summing the last two inequalities, multiplying by \(\gamk*\), and observing that
					\[
						\varphi_\star-\varphi(x^k)+\mytheta_{k+1}(\varphi(x^{k-1})-\varphi(x^k))
					=
						\mytheta_{k+1}P_{k-1}
						-
						(1+\mytheta_{k+1})P_k,
					\]
					we obtain
					\begin{align*}
						&
						\tfrac{1}{2}\|x^{k+1}-x^\star\|^2
						+
						\gamk*(1+\mytheta_{k+1})P_k
						+
						\tfrac{1-\varepsilon_{k+1}\rhok*}{2}\|x^k-x^{k+1}\|^2
					\\
					\leq{} &
						\tfrac{1}{2}\|x^k-x^\star\|^2
						+
						\mytheta_{k+1}\gamk*P_{k-1}
						+
						\rhok*\left(
							\tfrac{1-\gamk\Lk(2-\gamk\ck)}{2\varepsilon_{k+1}}
							-
							\mytheta_{k+1}(1-\gamk\Lk)
						\right)
						\|x^{k-1}-x^k\|^2.
					\numberthis\label{eq:PG:descentdelta}
					\end{align*}
					Selecting \(\mytheta_{k+1}=\rhok*\) and \(\varepsilon_{k+1}=\nicefrac{1}{2\rhok*}\) results in the claimed inequality.
					\qedhere
				\end{proof}
			\end{lemma}

			As detailed in the proof, the inequality in \cref{thm:PG:descent} is a special case of the more general \eqref{eq:PG:descentdelta} obtained by setting \(\mytheta_k=\rhok\) and \(\varepsilon_{k+1}=\nicefrac{1}{2\rhok*}\) for any \(k\).
			As we are about to see in the following \cref{thm:PG:convergence}, these choices strike a nice balance between simplifying inequality \eqref{eq:PG:descentdelta} and enabling large stepsizes: the result is a rather simple update rule for the stepsize that works very well in practice.
			It would be tempting to explore if more sophisticated tunings of these parameters could lead to further improvements, an aspect that we believe is worth investigating in the future, see \cref{rem:thetk}.

		\subsection{Convergence results}
			While \refadaPG* can be seen as a special case of the more general primal-dual \refadaPD*, the convergence results of \refadaPG* is obtained under fewer restrictions (cf. \cref{sec:PD:PG}).
			For this reason, we provide a dedicated proof for the adaptive proximal gradient algorithm.

			\begin{theorem}\label{thm:PG:convergence}%
				Suppose that \cref{ass:PG} holds, and consider the iterates generated by \refadaPG*.
				Then, for any \(x^\star\in\argmin\varphi\), with \(\rhok\) and \(P_k\) as in \eqref{eq:rhoP} and with \(\Uk(x^\star)\) defined as
				\[
					\Uk (x^\star) 				\coloneqq
					\tfrac12\|x^k-x^\star\|^2
					+
					\tfrac14\|x^k-x^{k-1}\|^2
					+
					\gamk\left(1+\rhok\right)P_{k-1},
				\]
				the following hold:
				\begin{enumerate}[widest*=4]
				\item \label{thm:PG:SD}%
					For all \(k\geq1\),
					\(
						\Uk*(x^\star)
					\leq
						\Uk (x^\star)
						-
						\bigl(
							\overbracket[0.5pt]{
								\tfrac{1}{4}
								-
								\rhok*^2\gamk\Lk(\gamk\ck-1)
							}^{\geq0}
						\bigr)
						\|x^{k-1}-x^k\|^2
						-
						\gamk\bigl(
							\overbracket[0.5pt]{
								\vphantom{\tfrac{1}{4}}
								{1+\rhok-\rhok*^2}
							}^{\geq0}
						\bigr)P_{k-1}
					\).
				\item \label{thm:PG:gamk}%
					The sequence \(\seq{x^k}\) is bounded and \(\gamk>\frac{1}{2L_{f,\V}}\) holds for every \(k\geq k_0\coloneqq\left\lfloor2\log_2\frac{1}{2\gamma_0L_{f,\V}}\right\rfloor_+\), where \(L_{f,\V}\) is a Lipschitz modulus for \(\nabla f\) on a compact convex set \(\V\) containing \(\seq{x^k}\).

				\item \label{thm:PG:x*}%
					The sequence \(\seq{x^k}\) converges to a solution.
					The claim remains true if the stepsizes are chosen in such a way that \eqref{eq:PG:gamk} holds with ``\(\leq\)'', as long as \(\seq{\gamk}\) is bounded away from zero.
					\item \label{thm:PG:sublin}%
					The following rate holds
					\[
						\min_{k=0,1,\ldots, K} (\varphi(x^k) - \min \varphi)
					\leq
						\frac{\U_1(x^\star)}{\sum_{k=1}^{K+1}\gamk}, 					\]
					which by assertion \itemref{thm:PG:gamk} implies a sublinear \(\mathcal O(\nicefrac1K)\) convergence rate on the best-so-far cost.
					\end{enumerate}
				\begin{proof}~
					\begin{proofitemize}
					\item \itemref{thm:PG:SD}~
						Expressed in terms of \(\U_k(x^\star)\), the inequality in \cref{thm:PG:descent} simplifies to the one in the statement. The update rule for \(\gamk\) implies that the coefficients indicated in the statement are greater than or equal to zero, and the claim follows.

					\item \itemref{thm:PG:gamk}~
						The proven inequality implies that \(\seq{\U_k(x_{}^{\star })}\) converges for any \(x^\star\in\argmin\varphi\).
						Taking into account the nonnegativity of the last two terms in the definition of \(\U_k(x^\star)\), it follows that \(\seq{x^k}\) is bounded.
						As such, \(\V\) and \(L_{f,\V}\) as in the statement exist.
							Observe that
							\begin{equation}\label{eq:relate:MM}
								\tfrac{\gamk}{2\sqrt{[\gamk\Lk(\gamk\ck-1)]_+}}
							=
								\tfrac12
								\sqrt{\tfrac{1}{\Lk\ck}\tfrac{\gamk\ck}{[\gamk\ck-1]_+}}
							\geq
								\tfrac{1}{2\sqrt{\Lk\ck}}
							=
								\tfrac{1}{2L_k},
							\end{equation}
							where the last identity owes to \cref{thm:PG:CL}, and therefore
							\begin{equation}\label{eq:MM}
								\gamk*
							\geq
								\min\set{\gamk\sqrt{1+\rhok},\tfrac{1}{2L_k}}
							\end{equation}
							holds for every \(k\).
												Suppose that \(\gamk\leq\nicefrac{1}{2L_{f,\V}}\) for \(k=0,\dots,K\).
							Because of \eqref{eq:MM} and \cref{thm:Lk<=Ck}, this implies that
							\(
								\gamk*
							=
								\gamk\sqrt{1+\rhok}
							\)
							for \(k=0,\dots,K-1\).
												Then, for these \(k\) one has that
							\(
								\rhok*^2
							=
								1+\rhok
							>
								1
							\),
							which inductively gives
							\(
								\rhok^2
							\geq
								2
							\)
							for all \(k=1,\dots,K\) (since \(\rho_0 \geq 1\)).
							We then have
							\begin{equation}\label{eq:PG:gam2k}
								\gamma_K^2
							=
								\rho_K^2\gamma_{K-1}^2
							\geq
								2\gamma_{K-1}^2
							\geq\dots\geq
								2^K\gamma_0^2, 							\end{equation}
							showing that in at most \(k_0\) iterations stepsize exceeds \(\nicefrac{1}{2L_{f,\V}}\).
							From this point, a trivial induction argument using \eqref{eq:MM} reveals that \(\gamk\geq\nicefrac{1}{2L_{f,\V}}\) for all \(k\geq k_0\).
					\item \itemref{thm:PG:x*}~
						We begin by observing that, as is apparent from its proof, assertion \itemref{thm:PG:SD} remains valid if the identity in \eqref{eq:PG:gamk} is replaced by an inequality ``\(\leq\)''.
						Either way, a telescoping argument yields that
						\begin{equation}\label{eq:gPto0}
							\gamk(1+\rhok-\rhok*^2)P_{k-1}\to0
						\quad\text{as}\quad
							k\to\infty.
						\end{equation}
						We proceed by intermediate claims.
						\begin{claims}[widest*=3]%
						\item \label{claim:infP}
							{\it \(\liminf_{k\to\infty}P_k=0\), and in particular \(\seq{x^k}\) admits a limit point \(x_\infty \in \argmin \varphi \).}

							Since \(\seq{x^k}\) is bounded and \(\varphi\) is lsc, it suffices to show that \(\liminf_{k\to\infty}P_k=0\).
							If \(\limsup_{k\to\infty}(1+\rhok-\rhok*^2)>0\), then \eqref{eq:gPto0} yields the claim, since \(\gamk\) is bounded away from zero.
							Alternatively, since \(1+\rhok-\rhok*^2\geq0\), necessarily \(1+\rhok-\rhok*^2\to0\), from which it easily follows that \(\liminf_{k\to\infty}\rhok>1\), and that therefore \(\gamk\to\infty\);
							in this case, because of the inequalities \(\gamk\left(1+\rhok\right)P_{k-1}\leq\U_k\leq\U_1\), this directly proves that \(P_k\to0\).

						\item\label{claim:thm:PG:optunique}%
							{\it \(x_\infty\) is the only limit point that belongs to \(\argmin\varphi\).}

							Suppose that \(x_\infty'\in\argmin\varphi\) is a limit point.
							Observe that
							\[
								\innprod{x^k}{x_\infty - x_\infty'}
							=
								\U_k(x_\infty') - \U_k(x_\infty) + \tfrac12\|x_\infty\|^2 - \tfrac12\|x_\infty'\|^2,
							\]
							and since \(\seq{\U_k(x^\star)}\) is convergent for all \(x^\star\in \argmin \varphi\), then so is \(\seq{\innprod{x^k}{x_\infty - x_\infty'}}\).
							Passing to the limit along the two converging subsequences thus yields
							\(
								\innprod{x_\infty}{x_\infty - x_\infty'}
							=
								\innprod{x_\infty'}{x_\infty - x_\infty'}
							\),
							which after rearranging results in \(\|x_\infty - x_\infty'\|^2=0\), establishing the claim.

						\item\label{claim:thm:PG:allopt}%
							{\it \(\seq{x^k}\) converges to a point in \(\argmin\varphi\).}

							Let \(x_\infty\) be the (unique) optimal limit point of \(\seq{x^k}\) as in \cref{claim:thm:PG:optunique}, and let
							\begin{equation}\label{eq:U}
								U\coloneqq\lim_{k\to\infty}\U_k(x_\infty),
							\end{equation}
							which exists since the sequence is monotonically decreasing by assertion \itemref{thm:PG:SD}.
							Because of the inequality \(\frac{1}{2}\|x^k-x_{\infty}\|^2\leq\U_k(x_{\infty})\), it suffices to show that \(U=0\).
							Observe that, for every \(k\),
							\begin{equation}\label{eq:rhomax}
								\rhok\leq\rho_{\rm max}
							\quad\text{where}\quad
							 	\rho_{\rm max}
							 \coloneqq
								\max\set{\tfrac{\gamma_0}{\gamma_{-1}}, \tfrac12(1+\sqrt{5})}
							>
								1,
							\end{equation}
							as it can be verified by induction.
							We consider two mutually exclusive scenarios.
							\begin{proofitemize}[leftmargin=*]
							\item
								Suppose that \(\seq{\gamk}\) is bounded, and consider a subsequence \(\seq{x^k}[k\in K]\) such that \(\lim_{K\ni k\to\infty}P_k=0\), which exists and converges to \(x_{\infty}\) by virtue of the previous claims.
								Extract \(K'\subseteq K\) such that \(\seq{\gamk*}[k\in K']\to\gamma\) for some \(\gamma>0\).
								Then, \(x^{k+1}=\prox_{\gamk*g}(x^k-\gamk*\nabla f(x^k))\to\prox_{\gamma g}(x_\infty-\gamma\nabla f(x_\infty))=x_\infty\) as \(K'\ni k\to\infty\) as well.
								This yields the sought identity
								\begin{align*}
									U
								={} &
									\lim_{k\to\infty}\U_k(x_{\infty})
								=
									\lim_{K'\ni k\to\infty}\U_{k+1}(x_{\infty})
								\\
								={} &
									\lim_{K'\ni k\to\infty}\left(
										\tfrac{1}2\|x^{k+1}-x_{\infty}\|^2
										+
										\tfrac{1}2\|x^{k+1}-x^k\|^2
										+
										\gamk*\left(1+\rhok*\right)P_k
									\right)
								\\
								={} &
									0,
								\end{align*}
								where the vanishing of the last term owes to boundedness of \(\seq{\gamk}\) and \(\seq{\rhok}\).

							\item
								Suppose instead that \(\seq{\gamk}\) is unbounded, and let \(\seq{\gamk}[k\in K']\) be an arbitrary subsequence such that \(\gamk\to\infty\) as \(K'\ni k\to\infty\).
								From the inequality
								\(
									\gamk P_{k-1} \leq \U_k(x_{\infty}) \leq \U_1(x_\infty) < \infty
								\),
								we conclude that
								\begin{equation}\label{eq:Pk:xk-}
									P_{k-1}\to0
								~\text{ and }~
									x^{k-1}\to x_\infty
								\text{ as }
									K'\ni k\to\infty
								\end{equation}
								(recall that \(\varphi\) is lsc, \(\seq{x^k}\) is bounded, and \(x_\infty\) is the \emph{unique} optimal limit point).
								Then, since by \eqref{eq:rhomax} \(\gamma_{k-1}\geq\nicefrac{\gamk}{\rho_{\rm max}}\to\infty\) as \(K'\ni k\to\infty\), by the same argument also \(x^{k-2}\to x_\infty\) as \(K'\ni k\to\infty\).
								With \(U\) as in \eqref{eq:U}, we thus have
								\begin{equation}\label{eq:U:gam:unbounded}
										U
								=
									\lim_{k\to\infty}\U_k(x_\infty)
								=
									\lim_{K'\ni k\to\infty}\U_{k-1}(x_\infty)
								=
									\lim_{K'\ni k\to\infty}\gamma_{k-1}(1+\rho_{k-1})P_{k-2}
								\end{equation}
								holding for \emph{any} set of indices \(K'\subseteq\N\) along which the stepsizes are divergent.

								We now construct a specific subsequence \(K\coloneqq\set{k_0,k_1,\dots}\) as follows: start with \(k_0=1\), and for \(i\geq0\) let
								\begin{equation}\label{eq:ki+}
									k_{i+1}=\min\set{k\geq k_i}[\gamk \geq \rho_{\rm max}\gamma_{k_i}].
								\end{equation}
								Then, \(\seq{\gamma_{k_i}}[i\in\N]\to\infty\), which as argued after \eqref{eq:Pk:xk-} implies that both \(\seq{x^{k_i-1}}[i\in\N]\) and \(\seq{x^{k_i-2}}[i\in\N]\) converge to \(x_\infty\).
								Observe that \(\rho_{k_i}>1\) holds for all \(i\in\N\) by minimality of \(k_i\).
								Notice further that
								\begin{equation}\label{eq:k-1}
									\rho_{k_i-1}\geq\rho_{\rm max}^{-1}
								\quad
									\forall i\in\N.
								\end{equation}
								To see why, suppose on the contrary that
								\(
									\rho_{k_i-1}<\rho_{\rm max}^{-1}
								\);
								then in particular \(\rho_{k_i-1}<1\), thereby ensuring by the previous observation that \(k_i-1\notin K\), and thus \(k_{i-1}\leq k_i-2\).
								This leads to the contradiction
								\[
									\rho_{\rm max}
									\gamma_{k_{i-1}}
								\overrel[\leq]{\eqref{eq:ki+}}
									\gamma_{k_i}
								\overrel[\leq]{\eqref{eq:rhomax}}
									\rho_{\rm max}
									\gamma_{k_i-1}
								=
									\rho_{\rm max}
									\rho_{k_i-1}
									\gamma_{k_i-2}
								<
									\gamma_{k_i-2}
								<
									\rho_{\rm max}
									\gamma_{k_{i-1}},
								\]
								where the last inequality follows either in case \(k_i-2=k_{i-1}\) (since \(\rho_{\rm max}>1\)) or when \(k_{i-1}<k_i-2<k_i\) (from minimality in the definition of \(k_i\)).
								This shows \eqref{eq:k-1}.

								Next, by the stepsize update at \cref{state:PG:gamk},
								\begin{align*}
									\tfrac{1}{\rho_{\rm max}}\gamma_{k_i-2}L_{k_i-2}
								\overrel[\leq]{\eqref{eq:k-1}}
									\gamma_{k_i-1}L_{k_i-2}
								\leq{} &
									\frac{\gamma_{k_i-2}L_{k_i-2}}{2\sqrt{\left[\gamma_{k_i-2}^2\ell_{k_i-2}c_{k_i-2}-\gamma_{k_i-2}\ell_{k_i-2}\right]_+}}
								\\
								\leq{} &
									\frac{\gamma_{k_i-2}L_{k_i-2}}{2\sqrt{\left[\gamma_{k_i-2}^2L_{k_i-2}^2-\gamma_{k_i-2}L_{k_i-2}\right]_+}}
								\\
								={} &
									\frac{1}{2\sqrt{\left[1-\frac{1}{\gamma_{k_i-2}L_{k_i-2}}\right]_+}},
								\end{align*}
								where the third inequality uses the fact that \(L_{k_i-2}^2 = \ell_{k_i-2}c_{k_i-2}\), and \(\ell_{k_i-2}\leq L_{k_i-2}\) by \cref{thm:PG:cLl}.
								This implies that \(\seq{\gamma_{k_i-2}L_{k_i-2}}[i\in\N]\) must be bounded.
								The subgradient characterization in \eqref{eq:PG:vk} and the Cauchy-Schwarz inequality then yield
								\begin{align*}
									\gamma_{k_i-1} P_{k_i-2}
								\leq{} &
									\rho_{k_i-1}
									\innprod{x^{k_i-2}-x_{\infty}}{
										x^{k_i-2}-x^{k_i-3}
										-
										\gamma_{k_i-2}(\nabla f(x^{k_i-2})-\nabla f(x^{k_i-3}))
									}
								\\
								\leq{} &
									\rho_{\rm max}
									\underbracket[0.5pt]{
										(1+\gamma_{k_i-2}L_{k_i-2})
									}_{\text{bounded}}
									\,
									\underbracket[0.5pt]{
										\|x^{k_i-2}-x_{\infty}\|
										\|x^{k_i-2}-x^{k_i-3}\|
									}_{\to0}
								\to
									0
								\quad
									\text{as }i\to\infty.
								\end{align*}
								Combined with the fact that \(\rho_{k_i-1} \leq \rho_{\rm max}\), it follows from \eqref{eq:U:gam:unbounded} that \(U=0\).
							\end{proofitemize}
												\end{claims}%

					\item \itemref{thm:PG:sublin}~
						A telescoping argument in the descent inequality of assertion \itemref{thm:PG:SD} yields
						\[
							\mathcal U_{K+1} + \sum_{k=1}^K\gamk\left(1+\rhok-\rhok*^2\right)P_{k-1} \leq \mathcal U_1
						\]
						for \(K\geq1\).
						By further lower bounding \(\mathcal U_{K+1}\) we obtain
						\begin{align*}
						\textstyle
							\gamma_{K+1}(1+\rho_{K+1})P_K
							+
							\sum_{k=1}^K\gamk\left(1+\rhok-\rhok*^2\right)P_{k-1}
							\leq \mathcal U_1.
						\end{align*}
						Considering the best-so-far quantity
						\(
							\min_{k=0,1,\ldots, K} P_k
						\), it follows that
						\begin{align*}
							\min_{k=0,1,\dots, K} P_k
						\leq{} &
							\frac{\mathcal U_1}{
								\gamma_{K+1}(1+\rho_{K+1})
								+
								\sum_{k=1}^K\gamk(1+\rhok-\rhok*^2)
							}
						\\
						={} &
							\frac{\mathcal U_1}{
								\gamma_{K+1}(1+\rho_{K+1})
								+
								\sum_{k=1}^K\gamk
								+
								\sum_{k=1}^K(\gamk\rhok-\gamk*\rhok*)
							}
						\\
						\numberthis\label{eq:O(1k)}
						={} &
							\frac{\mathcal U_1}{
								\gamma_1\rho_1+\sum_{k=1}^{K+1}\gamk
							}
						\leq
							\frac{\mathcal U_1}{
								(K+1)\gamma_{\rm min}
							},
						\end{align*}
						where \(\gamma_{\rm min}\coloneqq\min\set{\gamma_0,\frac{1}{2L_{f,\mathcal V}}}\) is a lower bound on the stepsize, as demonstrated in the proof of assertion \itemref{thm:PG:gamk}.
						\qedhere
					\end{proofitemize}
				\end{proof}
			\end{theorem}
			We remark that if the suggested initialization in \cref{sec:PG:init} is used, the initial stepsize satisfies \(\gamma_0 \geq \nicefrac{1}{L_{f,\V}}\) and consequently \eqref{eq:O(1k)} holds with \(\gamma_{\rm min}=\nicefrac{1}{2L_{f,\V}}\).
			The rate for the best-so-far cost thus simplifies to \(\min_{k=0,\dots,K}P_k\leq\frac{2\U _1L_{f,\V}}{K+1}\) in this case.
			More importantly, as a consequence of \cref{thm:PG:convergence} it can be easily seen that \cite[Alg. 1]{malitsky2020adaptive} can in fact cope with nonsmooth problems of the form \eqref{eq:PG}.
			This fact can be deduced from the inequality \eqref{eq:MM} in the proof of \cref{thm:PG:x*}, as formalized next.

			\begin{remark}[Comparison with {\cite[Alg. 1]{malitsky2020adaptive}}]\label{rem:PG:Malitsky}%
				With \(\gamk*\) as in \cref{state:PG:gamk} and \(L_k\) as in \eqref{eq:MML_k}, as shown in \eqref{eq:MM} one has that
				\begin{equation}\label{eq:adaPGvsMM}
					\gamk*
				\geq
					\min\set{
						\gamk\sqrt{1+\rhok},\,
						\tfrac{1}{2L_k}
					}
				\end{equation}
				holds for every \(k\). Note that the right-hand side corresponds to the stepsize update of \cite[Alg. 1]{malitsky2020adaptive}.
				As a consequence, in addition to accommodating proximal terms, \refadaPG* also comes with a less conservative stepsize update rule.
			\end{remark}

			Inequality \eqref{eq:adaPGvsMM} cannot be reiterated inductively, and in particular there is no guarantee that, iteration-wise, the \emph{entire} sequence of stepsizes produced by \refadaPG* is larger than that generated by \cite[Alg. 1]{malitsky2020adaptive}, even if the algorithms are started with same initial conditions.
			This is nevertheless enough to infer convergence of \cite[Alg. 1]{malitsky2020adaptive} applied to composite minimization problems as in \eqref{eq:PG}.

			\begin{corollary}[Proximal extension of {\cite[Alg. 1]{malitsky2020adaptive}}] \label{cor:MM}
				Suppose that \cref{ass:PG} holds.
				Then, {\algnamefont PGM} iterations \(x^{k+1}=\FBk*{x^k}\) with stepsize rule \eqref{eq:MML_k} converge to a solution of \eqref{eq:PG}.
				\begin{proof}
					The validity of \cref{thm:PG:x*} guarantees that the generated sequence remains bounded.
					As also observed in \cite[Thm. 1]{malitsky2020adaptive}, this guarantees that \(\gamk\geq\frac{1}{2L_{f,\V}}\) (up to possibly excluding initial iterates), where \(L_{f,\V}\) is a Lipschitz constant of \(\nabla f\) on a compact convex set \(\V\) that contains all the iterates.
					The proof follows by invoking \cref{thm:PG:x*} in light of \eqref{eq:adaPGvsMM}.
				\end{proof}
			\end{corollary}

			\begin{remark}[Alternative stepsize choices]\label{rem:thetk}%
				The update for \(\gamk\) in \refadaPG* is designed to ensure descent on the Lyapunov function \(\U_k\), cf. \cref{thm:PG:SD}. As commented after the proof of \cref{thm:PG:descent}, choices of the parameter \(\mytheta_k\) appearing in \eqref{eq:PG:descentdelta} other than \(\mytheta_k = \rhok\) lead to different update rules for the stepsize.
				For instance, retaining \(\varepsilon_k=\nicefrac{1}{2\rhok}\) but setting \(\mytheta_k=\pi\rhok\) for some \(\pi>0\) results in\footnote{%
					The update \eqref{eq:gamkpi} is obtained by expressing \eqref{eq:PG:descentdelta} in terms of
					\(
						\Uk^{\pi}(x^\star)
					\coloneqq
						\tfrac12\|x^k-x^\star\|^2
						+
						\tfrac14\|x^k-x^{k-1}\|^2
						+
						\gamk(1+\pi\rhok)P_{k-1},
					\)
					and enforcing descent as in \cref{thm:PG:SD}.%
				}%
				\begin{equation}\label{eq:gamkpi}
					\gamk*
				=
					\gamk\min\set{
						\sqrt{\tfrac1{\pi}+\tfrac{\gamk}{\gamma_{k-1}}},\,
						\frac{
							1
						}{
							2\sqrt{\left[\gamk\Lk(\gamk\ck-2+\pi) + 1 - \pi\right]_+}
						}
					},
				\end{equation}
				bringing about a trade-off between improving either term at the expense of the other.
				A similar concept is already pursued in \cite{malitsky2023adaptive}, which advances the update
				\[
					\gamk*
				=
					\gamk\min\set{
						\sqrt{\tfrac2{3}+\tfrac{\gamk}{\gamma_{k-1}}},\,
						\frac{
							1
						}{
							\sqrt{\left[2\gamk^2 L_k^2 - 1 \right]_+}
						}
					}
				\]
				(unrelated to \eqref{eq:gamkpi}, and obtained through different arguments) remarkably improving the worst-case rate coefficient, see \cite[Thm. 3 and \S 3.2]{malitsky2023adaptive}.
			\end{remark}

	\section{Adaptive three-term primal-dual methods}\label{sec:PD}
		In this section the idea of adaptively estimating the local geometry of \(f\) will be extended to composite problems of the form \eqref{eq:PD},	which we rewrite here for the reader's convenience
		\[
			\minimize_{x\in\R^n}\; \varphi(x) \coloneqq f(x) + g(x) + h(\linop x).
		\]
		Although the analysis in this setting is inevitably more complicated, the key idea of using the estimates \(\ck,\Lk\) in \eqref{eq:CL} remains the same.
		We will first propose an adaptive algorithm under the assumption that the norm of the linear operator \(\linop\) is known.
		This assumpion will then be lifted through a certain linesearch procedure.

		Problem \eqref{eq:PD} is studied under the following assumptions.

		\begin{assumption}[Requirements for the primal-dual setting]\label{ass:PD}%
			The following hold in problem \eqref{eq:PD}:
			\begin{enumeratass}
			\item
				\(\func{f}{\R^n}{\R}\) is convex and has locally Lipschitz continuous gradient.
			\item
				\(\func{h}{\R^m}{\Rinf}\) and \(\func{g}{\R^n}{\Rinf}\) are proper convex and lsc, and \(\func{\linop}{\R^n}{\R^m}\) is a linear mapping.
			\item
				A solution exists: \(\argmin \varphi \neq \emptyset\).
			\item \label{ass:CQ}%
				The problem is strictly feasible: there exists $x\in\relint \dom g$ such that $\linop x\in\relint \dom h$.
			\end{enumeratass}
		\end{assumption}

		A well-established approach for addressing \eqref{eq:PD} is to lift it into the primal-dual space and solve the associated convex-concave saddle point problem
		\[
			\minimize_{x\in\R^n} \maximize_{y\in\R^m} \SD(x,y)\coloneqq f(x)+g(x)+ \innprod{\linop x}y - \conj h(y).
		\]
		This lifting is the key to splitting the composed term \(h\circ\linop\).
		Moreover, in doing so the primal and the dual solutions can be obtained simultaneously.
		A pair \(z^\star=(x^\star,y^\star)\) will be referred to as a primal-dual solution if the following primal-dual optimality condition holds
		\begin{equation}\label{eq:PD:optInc}
				0 \in -\linop x^\star + \partial \conj h(y^\star),\quad
			0 \in \linop* y^\star + \nabla f(x^\star) + \partial g(x^\star).
		\end{equation}
		The set of all such pairs will be denoted by \(\mathcal S_\star\).
		Under the constraint qualification of \cref{ass:CQ}, the set of solutions for the dual problem is nonempty, and thus so is \(\mathcal S_\star\), and the duality gap is zero, see \cite[Cor. 31.2.1]{rockafellar1970convex} and \cite[Thm. 19.1]{bauschke2017convex}.
		Moreover, the pair $(x^\star,y^\star)$ is a primal-dual solution if and only if \(x^\star\) is a primal and \(y^\star\) is a dual solution.
		Since \eqref{eq:PD} is a convex problem, the primal-dual solution pairs are equivalently characterized by the saddle point inequality
		\begin{equation}\label{eq:saddleIneq}
			\SD(x^\star, y)
		\leq
			\SD(x^\star, y^\star)
		\leq
			\SD(x, y^\star)
		\quad
			\forall (x,y)\in\R^n\times\R^m.
		\end{equation}
		In our analysis we will measure deviation from \(\SD(x^\star,y^\star)\) along the primal and dual sequences using shorthand notations \(P_k=P_k(x^\star,y^\star)\) and \(Q_k=Q_k(x^\star,y^\star)\) defined as
		\begin{subequations}\label{subeq:PQ}%
			\begin{align}
			\label{eq:Pxy}
				P_k
			\coloneqq{} &
				\SD(x^k, y^\star)
				-
				\SD(x^\star, y^\star)
			=
				(f+g)(x^k)-(f+g)(x^\star) + \innprod{x^k-x^\star}{\linop* y^\star}
			\shortintertext{and}
				Q_k
			\coloneqq{} &
				\SD(x^\star,y^\star)
				-
				\SD(x^\star,y^k)
			=
				h^*(y^k) - h^*(y^\star) + \innprod{\linop x^\star}{y^\star-y^k},
			\end{align}
		\end{subequations}
		which are both positive due to the saddle point inequality \eqref{eq:saddleIneq}.

		\subsection{Algorithmic overview}\label{sec:PD:adaPDM}
			The proposed algorithm is presented in \refadaPD* (\cref{alg:PD}).
			It can be viewed as an adaptive variant of the algorithm proposed in \cite{condat2013primal,vu2013splitting}, which itself is an extension of the {\algnamefont PDHG} method \cite{chambolle2011first}.
			In comparison to the aforementioned algorithms with constant stepsizes, here a varying and potentially increasing stepsize rule is proposed based on the estimates \(\ck, \Lk\).
			Moreover, in {\algnamefont PDHG} the dual update combines two consequent primal updates as \(2 \linop x^k - \linop x^{k-1}\) (in our notation).
			In \cite[Alg. 3]{latafat2017asymmetric} it was shown that many primal-dual algorithms can be unified by modifying the dual update to use terms of the form  \(\theta\linop x^k+(1-\theta)\linop x^{k-1}\) followed by a correction step.
			While depending on the application this can lead to parallel implementations and potentially larger stepsizes compared to {\algnamefont PDHG} (with \(\theta = 2\)), the improvement in speed is limited by the use of global estimates.
			Differently from the aforementioned works, in \refadaPD* the mixing constant \(\theta\) is selected adaptively as \(1 + \nicefrac{\gamk*}{\gamk}\).
			When viewed as an extension of the proximal gradient method \refadaPG*, this idea appears natural.
			In fact, this was proposed in \cite{vladarean2021first} as a primal-dual extension of \cite{malitsky2020adaptive} where superior convergence rate compared to constant stepsize variants was observed.
			\refAdaPD* uses a different stepsize update rule and inherits the tighter estimates in \eqref{eq:CL} while permitting for the third nonsmooth term \(g\).
			Notice also that it is not symmetric with respect to the primal and dual variables and a different algorithm can be obtained by applying it to the dual problem.

			\begin{algorithm}[htb]
				\caption{Adaptive primal-dual method (\protect\adaPD)}
				\label{alg:PD}
				\begin{algorithmic}[1]
				\itemsep=3pt%

				\Require
					\begin{tabular}[t]{@{}l@{}}
						primal/dual (square inverse) stepsize ratio \(t>0\)
					\\
						stepsize parameters \(\sdcoeff>0\) and \(\const>1+\sdcoeff\)
						{\footnotesize (\eg, \(\sdcoeff=10^{-6}\), \(\const=1.2\))}
					\\
						initial primal-dual pair \((x^{-1},y^0)\in\R^n\times\R^m\)
						and
						stepsizes \(\gamma_{-1} \leq \gamma_0\leq\nicefrac{1}{2\const t\|\linop\|}\)
					\end{tabular}

				\Initialize
					\(
						x^0
					=
						\prox_{\gamma_0g}\bigl(
							x^{-1}-\gamma_0\nabla f(x^{-1})-\gamma_0\linop*y^0
						\bigr)
					\)
				\item[\algfont{Repeat for} \(k=0,1,\ldots\) until convergence]

				\State \label{state:PD:CL}%
					Set
					\(\delk=\gamk\Lk(\gamk\ck-1)\)
					and
					\(\xik=t^2\gamk^2\|\linop\|^2\),
					where \(\Lk\) and \(\ck\) are as in \eqref{eq:CL}

				\State \label{state:PD:gamk}%
					\begin{tabular}[t]{@{}l@{}}
						Define the stepsizes as
					\\
						\(
							\begin{cases}[@{}ll@{}]
								\gamk*
							={} &
								\min\set{
									\gamk\sqrt{1+\frac{\gamk}{\gamma_{k-1}}},
								~
									\frac{1}{2\const t\|\linop\|},
								~
									\gamk\sqrt{
										\frac{
											1-4\xik(1+\sdcoeff)^2
										}{
											2(1+\sdcoeff)
											\left(
												\sqrt{
													\delk^2
													+
													\xik(1-4\xik(1+\sdcoeff)^2)
												}
												+
												\delk
											\right)
										}
									}
								}
							\\
								\sigk*
							={} &
								t^2\gamk*\vphantom{\Big|_{\big|}}
							\end{cases}
						\)
					\end{tabular}

				\State \label{state:PD:y+}%
					\(
						y^{k+1}
					=
						\prox_{\sigk*\conj h}\left(
							y^k
							+
							\sigk*\left(
								\bigl(
									1+\tfrac{\gamk*}{\gamk}
								\bigr)\linop x^k
								-
								\tfrac{\gamk*}{\gamk}\linop x^{k-1}
							\right)
						\right)
					\)

				\State \label{state:PD:x+}%
					\(
						x^{k+1}
					=
						\prox_{\gamk*g}\bigl(
							x^k-\gamk*\nabla f(x^k)-\gamk*\linop*y^{k+1}
						\bigr)
					\)
				\end{algorithmic}
			\end{algorithm}

			\subsubsection{Initialization and practical considerations}\label{sec:PD:init}
				\refAdaPD* involves a few parameters, some of which require tuning in practice.
				The constant \(\sdcoeff\) is required to be strictly positive due to theoretical reasons, and in practice can be selected very close to zero (see \cref{rem:PG<PD:PQ} and the discussion thereafter).
				The constant \(\const\) affects the second and third arguments of the \(\min\) operator in \cref{state:PD:gamk}; too large a \(\const\) results in the second term limiting the stepsizes, while a value close to \(1+\sdcoeff\) can potentially result in the third term becoming too small during the next iterate.
				In practice, we suggest a value between 1.1 and 1.5 for this constant.
				Possibly a more critical parameter in play is \(t>0\), which denotes the ratio between the primal and the dual stepsizes.
				It is well known that this parameter can have a big impact on the performance of primal-dual methods.
				Even for the Condat-V\~u algorithm, where a simple stepsize condition \(\gamma \sigma \|\linop\|^2 \leq 1 - \tfrac{\gamma L_f}2\) is available, this parameter has to be tuned in general.
				Having an algorithm that adaptively selects the parameter \(t\) is an open research question.
				Finally, the initial primal stepsize can be chosen equal to the minimum between \(\nicefrac{1}{2\const t\|\linop\|}\) and a suitable estimate (obtained, \eg, as in \cref{sec:PG:init} for \refadaPG*).

			\subsubsection{Termination criteria}\label{sec:PDterm}%
				The primal and dual proximal updates lead to the following inclusions:
				\begin{subequations}\label{eq:opt:vk:PD}
					\begin{align}
					v_1^{k+1}
				\coloneqq{} &
					\tfrac{1}{\sigk*}(y^k-y^{k+1})
					+
					\tfrac{\gamk*}{\gamk}\left(\linop x^k-\linop x^{k-1}\right)+\left(\linop x^k-\linop x^{k+1}\right)
				\nonumber
				\\
				\in{} &
					\partial\conj h(y^{k+1})-\linop x^{k+1}
				\\
					v_2^{k+1}
				\coloneqq{} &
					\tfrac{1}{\gamk*}(x^k-x^{k+1})
					+
					\nabla f(x^{k+1})-\nabla f(x^k)
				\nonumber
				\\
				\in{} &
					\partial(g+f)(x^{k+1})+\linop*y^{k+1}.
				\end{align}
				\end{subequations}
				The sequence \(\seq{v^k}= \seq{v_1^k, v_2^k}\) is uniquely determined by the sequence \(\seq{z^k} = \seq{x^k,y^k}\) and can be computed efficiently without the need for additional gradient evaluations or matrix-vector products.
				Note that \(v^k \in Tz^k\), where
				\(
					T(x,y)
				\coloneqq
					\left(\partial h^*(y) - \linop x, \partial g(x) + \nabla f(x) + \linop* y \right)
				\)
				denotes the  operator associated with the primal-dual optimality conditions \eqref{eq:PD:optInc}. Therefore, the quantity
				\begin{equation}
					\|v^k\|\geq \dist(0, Tz^k) \label{eq:vk}
				\end{equation}
				serves as a measure of optimality and can be used as a termination criterion.

			\subsubsection{Comparison to the proximal gradient method}\label{sec:PD:PG}%
				When \(h\equiv0\) (and \(\linop=0\)), problem \eqref{eq:PD} reduces to problem \eqref{eq:PG} and \refadaPD* with \(\sdcoeff=0\) reduces to \refadaPG*.
				In fact, in this case apparently \(y^k\equiv 0\) and \(x^{k+1}=\FBk*{x^k}\), and since \(\xik\equiv0\) the stepsize update on \(\gamk\) reduces to
				\begin{align*}
					\gamk*
				={} &
					\gamk\min\set{
						\sqrt{1+\tfrac{\gamk}{\gamma_{k-1}}},
					~
						\tfrac{
							1
						}{
							\sqrt{
								2\bigl(
									|\delk|
									+
									\delk
								\bigr)^{\vphantom x}
							}
						}
					}
				\\
				={} &
					\gamk\min\set{
						\sqrt{1+\tfrac{\gamk}{\gamma_{k-1}}},
					~
						\tfrac{
							1
						}{
							2\sqrt{[\delk]_+}
						}
					},
				\end{align*}
				which is precisely the stepsize update rule in \refadaPG*.

				\begin{remark}\label{rem:PG<PD:PQ}%
					Although there is an apparent relation between \(P_k\) as in \eqref{eq:Pxy} and the one in \cref{thm:PG:descent}, while for \refadaPG* it was sufficient to show that \(\liminf_{k\to\infty } P_k = 0\) to infer existence of optimal limit points, in the primal-dual setting an argument through the cost function cannot be used.
					In fact, although
					\(\liminf_{k\to\infty} P_k = 0\) and \(\lim_{k\to\infty} Q_k = 0\) do hold even when \(\sdcoeff=0\) (cf. \cref{thm:PD:PQ}), and thus limit points \((\hat x, \hat y)\) of \(\seq{x^k,y^k}\) generated by either one of \cref{alg:PD,alg:PDls} exist that satisfy
					\[
						\SD(\hat x, y^\star)
					=
						\SD(x^\star,\hat y)
					=
						\SD(x^\star,y^\star)
						\quad
						\forall (x^\star,y^\star)\in \mathcal S_\star,
						\stepcounter{equation}
						\tag{\theequation}
						\label{eq:quasisaddle}
					\]
					we can only guarantee convergence to solutions by establishing sufficient descent in terms of the residual \(\|(x^{k+1},y^{k+1}) - (x^k,y^k)\|\) (enforced in the algorithms through the introduction of the parameter \(\sdcoeff>0\)).
				\end{remark}

				The fact that primal optimality of \(\hat x\) and/or dual optimality of \(\hat y\) cannot be inferred from \eqref{eq:quasisaddle} can be demonstrated with a simple counterexample.\footnote{%
					Example taken from \url{https://math.stackexchange.com/a/3039783/53739}.%
				}

				\begin{example}\label{ex:stackexchange}%
					Consider the Lagrangian \(\mathcal L(x,y)=x(1+y)\) of the problem
					\[
						\minimize_{x\in\R}\;x
					\quad
						\stt \;x=0.
					\]
					The unique saddle point of \(\mathcal L\) is \((x^\star,y^\star)=(0,-1)\), this being the only primal-dual solution of the problem.
					Nevertheless, any \(x\in\R\) minimizes \(\mathcal L(x,y^\star)\equiv0\) and similarly \(\mathcal L(x^\star,y)\equiv0\) is maximized by any \(y\in\R\);
					equivalently, \(P_k\) and \(Q_k\) as in \eqref{subeq:PQ} are identically zero independently of \(x^k\) and \(y^k\).
				\end{example}

				\Cref{ex:stackexchange} demonstrates that \(P_k\) and \(Q_k\) cannot, in general, be employed as optimality measures.
				Sufficient conditions involve the minimizer of \(\mathcal L({}\cdot{},y^\star)\) and/or the maximizer of \(\mathcal L(x^\star,{}\cdot{})\) being unique, which can be guaranteed under strict convexity assumptions on the primal and/or dual formulations.
				Whether \refadaPD* can accommodate \(\sdcoeff=0\) in the generality of \cref{ass:PD} remains an open question and is left for future work.

		\subsection{Linesearch variant without linear operator norm}\label{sec:PDls}
			The employment of the quantities \(\Lk\) and \(\ck\) in the adaptive stepsize strategies of \refadaPG* and \refadaPD* enables much tighter estimates of the local geometry of the problem, as opposed to adopting preset constants such as global Lipschitz moduli (when available).
			By the same principle, the norm of the linear operator \(\linop\) only offers a worst-case bound of the kind \(\|\linop*{}\cdot{}\|\leq\|\linop\|\|{}\cdot{}\|\), which can potentially be very loose on specific instances.
			To completely remove any dependency on global quantities, in this subsection we introduce \refadaPD+, a \emph{fully} adaptive primal-dual algorithm that replaces also the norm of the linear operator \(\linop\) with local estimates.
			As was also the case for the other two algorithms, \refadaPD+* follows the convention of indexing with \(k\) all the variables that depend on quantities defined up to iteration \(k\).
			This choice highlights the nested dependency of \(\gamk*\) and \(\eta_{k+1}\) at \cref{state:PD_eta:gamk}, which appears to be solvable only by means of a linesearch procedure (see also the discussion before \cref{thm:PD}).
			In account of this, we dub the method \emph{essentially} adaptive.
			The employment of a linesearch within an adaptive scheme is also pursued in \cite{latafat2023adabim}, where an extension of \refadaPG* is developed in the context of simple bilevel optimization, a setting that encompasses the primal-dual problem setting \eqref{eq:PD} as a special case.%
			 Backtracking linesearch has also been employed in combination with {\algnamefont PDHG} in various forms.
			In \cite{goldstein2013adaptive,goldstein2015adaptive} it is used to potentially increase the speed of convergence by balancing the primal and dual residuals.
			In \cite{malitsky2018first} an adaptive linesearch algorithm is presented for {\algnamefont PDHG} and the idea is extended to the composite form \eqref{eq:PD}, see \cite[Alg. 4]{malitsky2018first}.
			In addition to a different stepsize update rule, a major difference here is that the backtracks involved do not require evaluations of the gradient \(\nabla f\).
			We also remark that \refadaPD+* provides a practical way of initializing the stepsize, and that it reduces to \refadaPD* if \(\eta_k\) is taken as \(\|\linop\|\) for all \(k\).

			\begin{algorithm}[t]
				\caption{Adaptive primal-dual method with linesearch (\protect\adaPD+)}
				\label{alg:PDls}%
				\begin{algorithmic}[1]
				\itemsep=3pt%
				\renewcommand{\algorithmicindent}{0.3cm}%

				\Require
					\begin{tabular}[t]{@{}l@{}}
						primal/dual (square inverse) stepsize ratio \(t>0\)
					\\
						stepsize parameters \(\sdcoeff>0\) and \(\const>1+\sdcoeff\)
						{\footnotesize (\eg, \(\sdcoeff=10^{-6}\), \(\const= 1.2\))}%
					\\
						initial primal-dual pair \((x^{-1},y^0)\in\R^n\times\R^m\),~
						estimate \(\eta_0>0\) of \(\|\linop\|\),~
					\\
						stepsizes \(0<\gamma_{-1}\leq\gamma_0\leq\nicefrac{1}{2t\const\eta_0}\),~
						backtracking parameter \(r>1\)
						{\footnotesize (\eg, \(r=2\))}%
					\end{tabular}

				\Initialize
					\(x^0=\PDx[0]{x^{-1}}{y^0}\)

				\item[\algfont{Repeat for} \(k=0,1,\ldots\) until convergence]

				\State \label{state:PD_eta:CL}%
					Set \(\delk=\gamk\Lk(\gamk\ck-1)\) and
					\(
						\bar\xi_k
					=
						t^2\gamk^2\eta_k^2(1+\sdcoeff)^2
					\),
					with \(\Lk\) and \(\ck\) as in \eqref{eq:CL}%

				\State \label{state:PD_eta:hateta}%
					Choose an estimate \(0<\eta_{k+1}\leq\eta_k\) of \(\|\linop\|\) {\footnotesize (\eg, \(\eta_{k+1}=0.95\eta_k\))}

				\While{ \algfont{true} }
					\Comment{linesearch loop}

					\State \label{state:PD_eta:gamk}%
						\begin{tabular}[t]{@{}l@{}}
							Define the stepsizes as
						\\
							\(
								\begin{cases}[@{}l@{}l@{}]
									\gamk*
								={} &
									\min\set{\!
										\gamk\sqrt{1+\tfrac{\gamk}{\gamma_{k-1}}},
									\,
										\frac{1}{2\const t\eta_{k+1}},
									\,
										\gamk
										\sqrt{
											\frac{
												1-4\bar\xi_k
											}{
												2(1+\sdcoeff)
												\left(
													\sqrt{
														\delk^2
														+
														(t\eta_{k+1}\gamk)^2(1-4\bar\xi_k)
													}
													+
													\delk
												\right)
												\vphantom{X^{X^X}}
											}
										}
										\!
									}
								\\
									\sigk*
								={} &
									t^2\gamk*\vphantom{\Big|_{\big|}}
								\end{cases}
							\)
						\end{tabular}

					\State \label{state:PD_eta:y+}%
						\(
							y^{k+1}
						=
							\prox_{\sigk*\conj h}\left(
								y^k
								+
								\sigk*\left(
									\bigl(
										1+\tfrac{\gamk*}{\gamk}
									\bigr)\linop x^k
									-
									\tfrac{\gamk*}{\gamk}\linop x^{k-1}
								\right)
							\right)
						\)

					\State \label{state:PD_eta:break}%
						\algfont{if }
						\(
							\eta_{k+1}
						\geq
							\frac{\|\linop*(y^{k+1}-y^k)\|}{\|y^{k+1}-y^k\|}
						\)
						\algfont{ then break, ~else }
						\(\eta_{k+1}\gets r\eta_{k+1}\)

				\EndWhile

				\State \label{state:PD_eta:x+}%
					\(x^{k+1}=\PDx*{x^k}{y^{k+1}}\)
				\end{algorithmic}
			\end{algorithm}

			\subsubsection{Practical considerations and termination criteria}\label{sec:PD+term}%
				The same initialization and parameter values described for \refadaPD* can be used in \refadaPD+* with the difference that the initial estimate \(\eta_0\) replaces \(\|\linop\|\).
				The estimate \(\eta_0\) can be computed, for instance, by evaluating the Frobenius norm of the linear operator as an underestimation for it.
				Subsequent values of \(\eta_{k+1}\) can be initialized based on the previously accepted value \(\eta_k\); consistently with what suggested in \cref{state:PD_eta:hateta}, a multiple in the range \([0.9, 1)\) is recommended.
				Regarding the termination criterion, the same strategy discussed in \cref{sec:PD:adaPDM} for \refadaPD* can be employed.

		\subsection{Convergence results}\label{sec:convergence:PD}
			The convergence analysis will once again revolve around showing descent on a suitable merit function, this time defined as \(\Uk=\Uk(x^\star,y^\star)\) given by
			\begin{equation}\label{eq:PD:U}
				\Uk
			\coloneqq
				\tfrac{1}{2}\|x^k-x^\star\|^2
				+
				\tfrac{1-4\xik(1+\sdcoeff)}{4}
				\|x^k-x^{k-1}\|^2
				+
				\tfrac{1}{2t^2}\|y^k-y^\star\|^2
				+
				\gamk(1+\rhok)P_{k-1},
			\end{equation}
			where \((x^\star,y^\star)\in\mathcal S_\star\) is any primal-dual optimal pair.
			The next theorem allows us to study the convergence of \cref{alg:PD,alg:PDls} under a unified analysis.
			It should be noted that, unless \(\|\linop\|\) is known and \(\eta_{k+1}\) is chosen greater to or equal than that quantity (as it happens in \refadaPD*), the following theorem does not furnish an implementable stepsize update rule, owing to the implicit dependency between the stepsize \(\gamk*\) and the norm estimate \(\eta_{k+1}\).
			The linesearch prescribed by \refadaPD+* circumvents this issue.

			\begin{theorem}\label{thm:PD}%
				Suppose that \cref{ass:PD} holds, and let \(t>0\), \(\sdcoeff\geq0\), and \(\const>1+\sdcoeff\) be fixed.
				Consider a sequence \(\seq{x^k,y^k}[k\geq1]\) generated by
				\begin{equation}\label{eq:PD:update}
					\renewcommand{\arraystretch}{1.4}
					\begin{cases}[rl]
						y^{k+1}
					={} &
						\prox_{\sigk*\conj h}\left(
							y^k
							+
							\sigk*\left(
								\bigl(
									1+\tfrac{\gamk*}{\gamk}
								\bigr)\linop x^k
								-
								\tfrac{\gamk*}{\gamk}\linop x^{k-1}
							\right)
						\right)
					\\
						x^{k+1}
					={} &
						\PDx*{x^k}{y^{k+1}},
					\end{cases}
				\end{equation}
				starting from a triplet \((x^{-1},x^0,y^0)\in\R^n\times\R^n\times\R^m\) and with initial primal stepsizes \(\gamma_0\geq \gamma_{-1}>0\).
				Denote \(\xik\coloneqq t^2\eta_k^2\gamk^2\) and \(\delk\coloneqq\gamk\Lk(\gamk\ck-1)\) with \(\Lk\) and \(\ck\) as in \eqref{eq:CL}, \(\eta_0\leq\frac{1}{2\const t\gamma_0}\) and \(\eta_{k+1}\) any such that \(\frac{\|\linop*(y^k-y^{k+1})\|}{\|y^k-y^{k+1}\|}\leq\eta_{k+1}\leq\eta_{\rm max}\) for some \(\eta_{\rm max}<\infty\), \(k\in\N\).
				Suppose that the sequences of stepsizes comply with the rules
				\[
				\textstyle
					\gamk*
				=
					\min\set{\!
						\gamk\sqrt{1+\tfrac{\gamk}{\gamma_{k-1}}},
					\,
						\frac{1}{2\const t\eta_{k+1}},
					\,
						\gamk
						\sqrt{
							\frac{
								1-4\xik(1+\sdcoeff)^2
							}{
								2(1+\sdcoeff)
								\left(
									\sqrt{\delk^2+(t\eta_{k+1}\gamk)^2(1-4\xik(1+\sdcoeff)^2)\,}
									+
									\delk
								\right)
							}
						}
					\!}
				\]
				and \(\sigk=t^2\gamk\), \(k\in\N\).
				Then, for any primal-dual solution \((x^\star,y^\star)\) of \eqref{eq:PD} and with \(\Uk\) as in \eqref{eq:PD:U}, the following hold:
				\begin{enumerate}[topsep=0pt]
				\item\label{thm:PD:SD}%
					\(
						\Uk*
					\leq
						\Uk
						-
						\tfrac{\sdcoeff}{2t^2(1+\sdcoeff)}
						\|y^k-y^{k-1}\|^2
						-
						\underbracket*[0.5pt]{
							\tfrac{1-4\xik-4\rhok*^2(\delk+\xik*)}{4}
						}_{
							\geq\nicefrac{\sdcoeff}{4(1+\sdcoeff)}
						}
						\|x^k-x^{k+1}\|^2
						-
						\gamk(
							\overbracket[0.5pt]{
								\vphantom{\tfrac14}
								1+\rhok-\rhok*^2
							}^{\geq0}
						)P_{k-1}
						-
						\gamk*Q_{k+1}
					\).

				\item \label{thm:PD:gamk}%
					The sequence \(\seq{x^k,y^k}\) is bounded, and \(\gamk \geq \hat \gamma >0\) for all \(k\geq1\) (see \eqref{eq:PDgamkLB} for the value of \(\hat\gamma\)).

				\item \label{thm:PD:PQ}%
					\(\liminf_{k\to\infty}P_k=\lim_{k\to\infty}Q_k=0\).

				\item \label{thm:PD:z*}%
					If \(\sdcoeff>0\), the sequence \(\seq{x^k,y^k}\) converges to a primal-dual solution.%
				\end{enumerate}
			\end{theorem}

			We also remark that, thanks to the descent inequality in \cref{thm:PD:SD}, the same telescoping arguments as in \cite{chambolle2011first} can be used to show an \(\mathcal O(\nicefrac1k)\) convergence rate in terms of the partial primal-dual gap function introduced in \cite{chambolle2011first}.
			The sequential convergence results stated next for \refadaPD* follow from the more general \hyperref[thm:PD:z*]{Theorem }\ref{thm:PD:z*}, specialized to \(\eta_k\equiv\eta_{\rm max}=\|\linop\|\).
			For \refadaPD+* the assertions follow from \cref{thm:PD} this time with \(\eta_{\rm max}=\max\set{\eta_0,r\|\linop\|}\), as this furnishes an upper bound on \(\seq{\eta_k}\) (this fact simply follows by observing that whenever \(\eta_{k+1}\geq\|A\|\) one has that the backtracking procedure terminates at \cref{state:PD_eta:break}, and that \(\eta_{k+1}\) is otherwise increased by a factor \(r\) at any failed attempt).

			\begin{theorem}[convergence of \cref{alg:PD,alg:PDls}]%
				Under \cref{ass:PD}, all the assertions in \crefrange{thm:PD:SD}{thm:PD:z*} remain valid for both of the sequences generated by \cref{alg:PD,alg:PDls}.
			\end{theorem}

	\section{Numerical simulations}\label{sec:simulations}
		In this section the performance of the proposed algorithms is evaluated through a series of simulations on standard problems on both synthetic data as well as datasets from the LIBSVM library \cite{chang2011libsvm}.
		All the algorithms are implemented in the Julia programming language and are available online.\footnote{%
			\url{https://github.com/pylat/adaptive-proximal-algorithms}%
		}

			\begin{table}[h]
				\centering%
				\begin{tabular}{@{}>{\algnamefont}ll@{}}
					PGM       & Proximal gradient with constant stepsize {$\nicefrac{1}{L_f}$}%
				\\
					PGM-ls    & Proximal gradient with backtracking
				\\
					Nesterov  & Nesterov's acceleration with constant stepsize $\nicefrac{1}{L_f}$ \cite[\S10.7]{beck2017first}%
				\\
					aGRAAL    & The golden ratio algorithm \cite{malitsky2020golden}%
				\\
					PDHG      & The algorithm of \cite{chambolle2011first}%
				\\
					CV        & The algorithm of Condat and V\~u, proposed in \cite{condat2013primal,vu2013splitting}%
				\\
					MP-ls     & The linesearch method of Malitsky and Pock \cite[Alg. 4]{malitsky2018first}%
				\\
					\adaPG    & \cref{alg:PG}%
				\\
					\adaPG-MM & Proximal extension of Malitsky and Mishchenko \cite[Alg. 1]{malitsky2020adaptive}\footnotemark
				\\
					\adaPD    & \cref{alg:PD}%
				\\
					\adaPD+   & \cref{alg:PDls}%
				\end{tabular}
				\caption{Algorithms compared against in the numerical simulations (when applicable).}%
			\end{table}
		\footnotetext{%
			See \cref{cor:MM}.%
		}%

		The backtracking procedure in {\algnamefont PGM-ls} is meant in the sense of \cite[\S10.4.2]{beck2017first}, (see \cite[LS1]{salzo2017variable} and \cite[Alg. 3]{demarchi2022proximal} for the locally Lipschitz smooth case), without enforcing monotonic decrease on the stepsize sequence.
		Specifically, the initial guess for \(\gamk*\) is warm-started as \(r\gamk\), where \(\gamk\) is the accepted value in the previous iteration and \(r\geq1\) is a scaling factor.
		In each simulation the best outcome for {\algnamefont PGM-ls} among the choices \(r\in\set{1, 1.5, 2}\) is reported.

		We also considered two linesearch variants of Nesterov's accelerated method.
		The first variant is described in \cite[Eq. 4.9]{nesterov2013gradient} and allows the stepsize to increase (by warm starting the backtracks with a multiple of the last accepted stepsize).
		The second variant is described in \cite[\S10.7]{beck2017first} and initializes the backtracking procedure with the last accepted stepsize.
		We tested the first algorithm with scaling factors of 1.5 and 2.
		Both variants performed worse than the constant stepsize regime, and we have therefore omitted them from the plots.
		The degradation in performance is likely due to the extra cost of evaluating the backtracking condition.

		\subsection{Adaptive proximal gradient}
			We compare the performance of \refadaPG{} on three practical problems that can be cast as \eqref{eq:PG}.
			In the figures, the distance of the cost from the minimum is plotted against the number of calls to linear operations, which in all problems correspond to the most costly operation.
			This accounts for all gradient evaluations, and all additional cost evaluations that linesearch methods incur.

\begin{figure}[htbp]
				\begin{subfigure}{\linewidth}
					\includetikz[width=\linewidth]{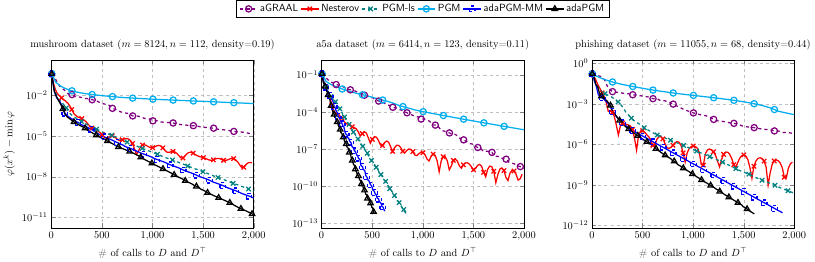}%
					\caption{Logistic regression problem of \hyperref[sec:logreg]{Section} \ref{sec:logreg}}
					\label{fig:logreg}%
				\end{subfigure}
				\begin{subfigure}{\linewidth}
					\includetikz[width=\linewidth]{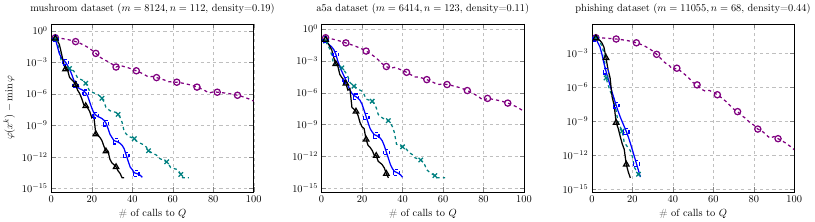}%
					\caption{Cubic regularization problem of \hyperref[sec:cubic]{Section} \ref{sec:cubic}}
				\end{subfigure}
				\begin{subfigure}{\linewidth}
					\includetikz[width=\linewidth]{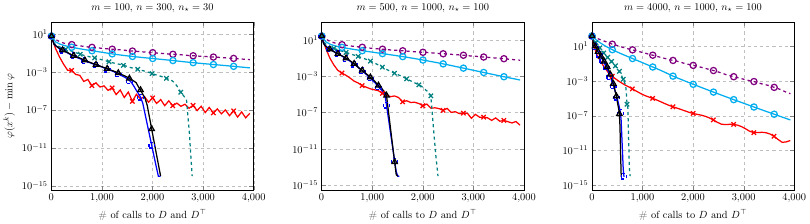}%
					\caption{Lasso problem of \hyperref[sec:lasso]{Section} \ref{sec:lasso}}
				\end{subfigure}
				\caption{%
					Simulations for problem \eqref{eq:PG}.
				}%
				\label{fig:PG}%
			\end{figure}

			In order to have a fair comparison between different methods, for globally Lipschitz-smooth problems we initialized the stepsize for all the methods with \(\nicefrac{1}{L_f}\) where \(L_f\) denotes the Lipschitz constant of \(\nabla f\).
			For cubic regularization we perturbed the initial point with a random vector and used the inverse of the estimate as in \eqref{eq:MML_k} using these two points.
			We note that \(r=1.5\) usually worked best in backtracking linesearch variants; for other parameters we used standard choices.
			For {\algnamefont aGRAAL}, we set the algorithm parameters as suggested in \cite[\S 5]{malitsky2020golden}.
			The inital point \(x^0\) was set to zero for all algorithms.
			The results of all the simulations are reported in \cref{fig:PG}.

			\subsubsection{Logistic regression}\label{sec:logreg}
				We consider an \(\ell_1\)-regularized logistic regression problem
				\[
					\minimize_{x\in \R^{n+1}} \; -\tfrac1m
					\textstyle
					\sum_{i=1}^m
					\Bigl[
						y_i \log\bigl(s_i(x)\bigr)
						+
						(1-y_i)\log\bigl(1-s_i(x)\bigr)
					\Bigr]
					+ \lambda \|x\|_1,
				\]
				where \(\lambda>0\) is the regularization parameter, \(m,n\) are the number of samples and features,
				\(
					s_i(x)
				=
					(1+ \exp(-\trans{d_i}x))^{-1}
				\)
				are logistic sigmoid functions with \(d_i\in\R^{n+1}\) denoting the \(i\)-th data sample (up to absorbing the bias terms), and \(y_i\in\set{0,1}\) is the associated label.
				In \cref{fig:logreg} the algorithms are compared in terms of the total number of calls to \(D\) and its adjoint, where \(D\) denotes the data matrix constructed by stacking the vectors \(d_i\).
				In all the simulations \(\lambda = 0.01\) was used.

			\subsubsection{Cubic regularization}\label{sec:cubic}
				The subproblem in the cubic Newton method \cite{nesterov2006cubic} involves solving
				\[
				    \minimize_{x\in \R^{n+1}} \; \tfrac12\innprod{x}{Qx} + \innprod{x}{q} + \tfrac{M}{6}\|x\|^3,
				\]
				where \(Q\in\R^{(n+1)\times (n+1)}\), \(q\in\R^{n+1}\), and \(M>0\) is some regularization parameter.
				In the simulations, the Hessian \(Q\) and the gradient \(q\) are generated for the logistic loss problem evaluated at zero on the mushroom, a5a, and phishing datasets.
				Moreover, \(M= 1\) is used in the plots, though we remark that the behavior of the algorithms for different values is similar.

			\subsubsection{Regularized least squares}\label{sec:lasso}
				Consider the lasso problem
				\[
					\minimize_{x\in\R^n}\tfrac12\|Dx - b\|^2 + \|x\|_1,
				\]
				where matrix \(D\in \R^{m\times n}\) and vector \(b\in \R^n\) are generated based on the procedure described in \cite[\S6]{nesterov2013gradient}, and \(n_\star\) denotes the number of nonzero elements of the solution.
				In all simulations, the parameter controlling the magnitude of the primal solution was set equal to one (\(\rho = 1\) in the notation of \cite[\S6]{nesterov2013gradient}), but the general behavior of the algorithms is similar for alternative values.

		\subsection{Adaptive primal-dual algorithms}
			We now compare the performance of \refadaPD{} and its operator norm--free extension \refadaPD+ on problems fitting into formulation \eqref{eq:PD}.

			\subsubsection*{Optimality criterion}
				The optimality criterion defined in \cref{sec:PDterm} will be used in all the simulations for comparison between algorithms.
				It consists of computing the norm of the  sequence \(v^k\) uniquely determined by the generated sequence for which  \(v^k\in Tz^k\) holds.
				This quantity is provided in \eqref{eq:opt:vk:PD} for \refadaPD* and \refadaPD+*, as well as the Condat-V\~u and {\algnamefont PDHG} algorithms (both corresponding to updates in \cref{state:PD:y+,state:PD:x+} with constant stepsize).
				In the case of {\algnamefont MP-ls} \cite[Alg. 4]{malitsky2018first}, it is similarly obtained according to the optimality conditions of the proximal updates at steps 1 and 2.a therein (the extra evaluations needed are not counted in the plots).

			\subsubsection*{Stepsize selection}
				The stepsize parameters for \refadaPD* and \refadaPD+* were configured as follows:
				\(\sdcoeff = 10^{-8}\),
				\(\const=1.2\),
				and
				\(\gamma_0=\gamma_{-1} = 1 / (2 \const t \eta_0)\).
				In the presented simulations, for comparison purposes we used \(\eta_0=\|\linop\|\), noting that \refadaPD+* displays low sensitivity to the initial choice of \(\eta_0\).
				Consenquently, there is little need for fine-tuning this parameter and any reasonable estimate results in a very similar trajectory.

				We observed that the stepsize ratio \(t\) requires some tuning.
				In our preliminary simulations for \cref{alg:PD,alg:PDls}, on the problems presented in here, we ran a grid search for \(t\in [0.01, 100]\) and observed that often \(t=1\) performed well enough on the presented examples.
				For the Condat-V\~u and {\algnamefont PDHG} algorithms, we used the heuristic stepsize selection rule suggested in \cite[Eq. (3.53)]{latafat2020distributed}.
				In the case of {\algnamefont MP-ls}, in  \cite[\S5]{malitsky2018first} it was suggested to set \(t\) equal to the ratio obtained by tuning the constant stepsize regime of the Condat-V\~u algorithm.
				In our experiments, in addition to this heuristic the performance of {\algnamefont MP-ls} was fine-tuned through a grid search in the range \(t\in[0.01, 100]\).

			\subsubsection{Dual support vector machine problem}\label{sec:dsvm}%
				The support vector machine (SVM) problem consists of solving
				\[
					\minimize_{x\in\R^n}\;
					C\sum_{n=1}^N \max\set{ 0, 1- a_i (\trans{d_i}x + x_0)} + \tfrac12\|x\|^2,
				\]
				where \((d_i,a_i)\in\R^n\times\R\) represents the \(i\)-th data pair, and \(C>0\) is some positive constant.
				\begin{figure}[tb]
					\includetikz[width=\linewidth]{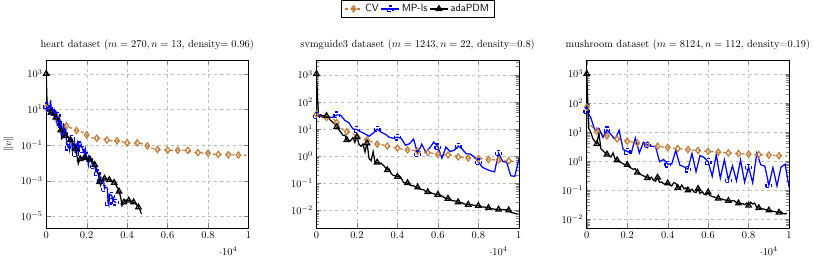}
					\includetikz[width=\linewidth]{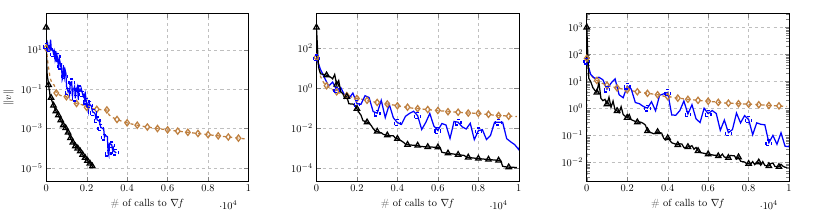}%
					\caption[]{%
						Simulations for the dual SVM problem \eqref{eq:dsvm}.
						First row: \(C = 1\), second row: \(C = 0.1\).
						The \(x\)-axis reports gradient evaluations, which is the most expensive operation (since \(\protect\linop\in \R^{1\times N}\), calls to \(\protect\linop\) and \(\protect\linop*\) are negligible).
						As explained, in this case {\protect\refadaPD+*} is indistinguishable from {\protect\refadaPD*} and thus omitted.
						\refAdaPD* and {\algnamefont MP-ls} are tuned for best performance by a grid search for \(t\in[0.01, 100]\).%
					}%
					\label{fig:dsvm}%
				\end{figure}
				Computationally, a popular approach is to instead consider the dual SVM problem \cite[\S12.2.1]{hastie2001elements}
				\begin{align*}
					\minimize_{\alpha_1,\ldots,\alpha_N}
				\;&
					\textstyle
					\tfrac{1}{2}\|\sum_{i=1}^N \alpha_i a_id_i\|^2-\sum_{i=1}^N \alpha_i
				\\
				\numberthis\label{eq:dsvm}
					\stt
				\;&
					0\leq\alpha_i\leq C,\quad i=1,\ldots,N
				\\
				&
					\textstyle\sum_{i=1}^N \alpha_i a_i = 0.
				\end{align*}
				In the simulations we used \(C\in\set{1, 0.1}\).
				The problem is cast in the standard form \eqref{eq:PD} by letting \(f\) represent the quadratic cost, \(g\) the indicator of the box constraints, \(h=\indicator_{\set{0}}\), and \(\linop = [a_1, \dots, a_N]\in \R^{1\times N}\).
				Notice that the dual vectors \(y^k\) are scalars in this case, and in particular one has \(\|\linop*(y^{k+1}-y^k)\|=\|\linop\|\|y^{k+1}-y^k\|\) for any \(k\), where \(\|\linop\|\) is a vector norm which is negligible to compute.
				Apparently, in this case (up to dicarding inital iterates) the linesearch variant \refadaPD+* produces the same iterates of \refadaPD* with no computational advantage, and is thus omitted from the plots.
				The results of this simulation are reported in \cref{fig:dsvm}.

			\subsubsection{Least absolute deviation regression and square-root lasso}\label{sec:MDreg}
				As a final application we consider regularized regression, consisting of
				\begin{equation}\label{eq:medreg}
					\minimize_{x\in\R^n}\|Dx-b\|_p + \lambda\|x\|_1.
				\end{equation}
				\begin{figure}[t]
					\includetikz[width=\linewidth]{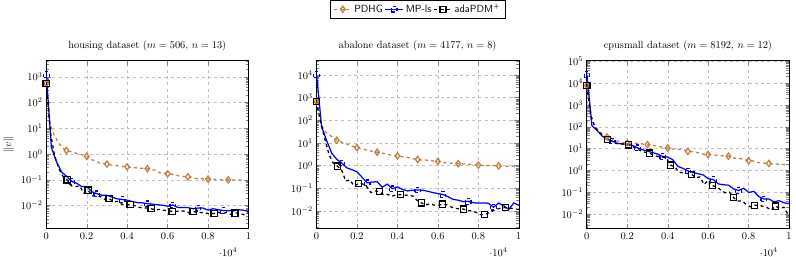}
					\includetikz[width=\linewidth]{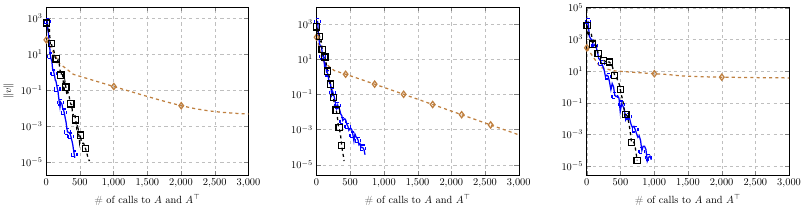}%
					\caption[]{%
						Simulations for problem \eqref{eq:medreg} of \cref{sec:MDreg} with \(\lambda = 10\).
						First row: regularized least-absolute deviation (\(p=1\));
						second row: square-root lasso (\(p=2\)).
						{\algnamefont AdaPDM\textsuperscript+} and {\algnamefont MP-ls} are tuned for best performance by a grid search for \(t\in[0.01, 100]\).
						Having \(f = 0\), \adaPD{} reduces to {\algnamefont PDHG} with worse (constant) stepsizes and is thus omitted from the comparisons.%
					}%
					\label{fig:medreg}
				\end{figure}
				Whenever \(p=2\), this problem is referred to as square-root lasso, and for \(p=1\) it is known as the least absolute deviation (LAD) regression.
				Both problems are cast as \eqref{eq:PD} by letting \(f\equiv 0\), \(g = \|\cdot\|_1\), \(h = \|{}\cdot{} - b\|_p\), and \(\linop = D\).
				For the regression data \(D\in \R^{m\times n}\) and \(b\in\R^m\), three different datasets from the LIBSVM library were used.
				The results are reported in \cref{fig:medreg}.
				It is apparent that the linesearch variant \refadaPD+* (as well as {\algnamefont MP-ls}, for similar reasons) excels by adaptively employing tighter bounds of the linear operator norm along one-dimensional subspaces.
				We ran the simulations for \(\lambda \in\set{0.1,1,10}\), but only report them for \(\lambda = 10\) remarking that the behavior is similar for the other values.

		\subsection{Observations about the stepsize sequence}\label{sec:observations}
			\Cref{thm:PG:gamk} certifies that the stepsizes \(\gamk\) produced by \refadaPG* do not drop below \(\nicefrac{1}{2L_{f,\V}}\), where \(L_{f,\V}\) is a Lipschitz modulus for \(\nabla f\) on a convex and bounded set \(\V\) that contains all the iterates.
			A refined analysis can nevertheless reveal some interesting additional details that, informally, confirm the pattern exhibited by the sequence \(\seq{\gamk}\) in all our experiments.

			\begin{figure}[ht]
				\centering
				\includetikz[width=.66\linewidth]{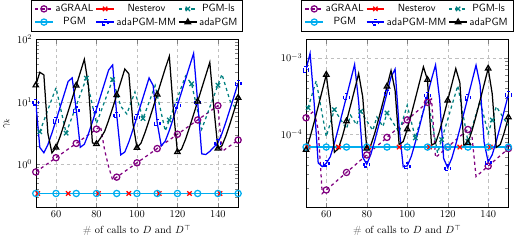}
				\includetikz[width=.66\linewidth]{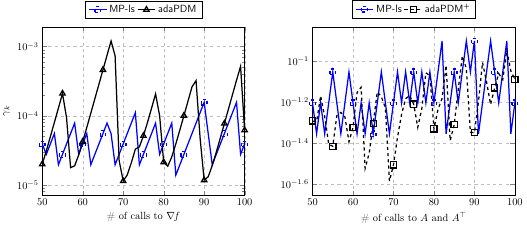}%
				\caption[]{%
					Demonstrative plot of stepsize magnitudes in windows of 50 and 100 iterations extracted from the simulations in this section.
					First row, left: logistic regression, mushroom dataset;
					right: lasso problem with \(m=500\), \(n=1000\), \(n_\star = 100\).
					As commented after \cref{rem:PG:Malitsky}, despite the fact that the stepsize update rule of \adaPG{} is less conservative than that of {\algnamefont\adaPG-MM}, the comparison does not carry over iterationwise.
					Second row, left: dual SVM, mushroom dataset, \(C=0.1\);
					right: square-root lasso (housing dataset, \(\lambda=10\)).
					Primal-dual algorithms are compared on problems where the ratio \(t^2=\sigk/\gamk\) coincides, so that the plots are representative also for dual stepsizes \(\sigk\).
				}%
				\label{fig:gam}%
			\end{figure}

			Whenever a stepsize is below a certain threshold, a first phase is triggered where the stepsize strictly increases until the threshold is surmounted.
			As ensured by \cref{lem:enlarged:ck}, \(\ck\) is upper bounded by a local Lipschitz modulus \(L_{f,\overline{\V}}\) for \(\nabla f\) on some enlargement \(\overline{\V}\) of the set \(\V\) containing the iterates.
			Because of the update rule at \cref{state:PG:gamk}, if \(\gamk\leq\nicefrac{1}{L_{f,\overline{\V}}}\) (hence \(\gamk\ck\leq1\)) then \(\gamk*>\gamk\) (in fact it increases at least linearly, cf. \eqref{eq:PG:gam2k}).
			This is followed by a second phase of unspecified length, during which the stepsizes remain strictly greater than \(\nicefrac{1}{L_{f,\overline{\V}}}\).
			Should a stepsize again drop below it (yet surely not below \(\nicefrac{1}{2L_{f,\V}}\)), the cycle restarts from phase one.
			Such an oscillatory behavior is confirmed in all our simulations (see \cref{fig:gam} for particular instances) and could also explain the adoption of much larger stepsizes compared to the standard constant stepsize regime.
			Whether more rigorous arguments could be made in support of this empirical evidence is an interesting open question left for future work.

	\section{Conclusions}\label{sec:conclusions}
		In this paper we studied convex composite problems involving the sum of a locally Lipschitz differentiable term and two nonsmooth terms, one of which composed with a linear operator.
		The primal-dual algorithm \refadaPD* was proposed that updates the stepsizes adaptively using already available information without any backtracking procedures.
		The stepsize update is based on a novel rule that combines estimates of local cocoercivity and Lipschitz continuity of the gradient of the differentiable function.
		When the linear operator is zero we obtain \refadaPG* that not only extends the adaptive gradient descent algorithm of \cite{malitsky2020adaptive}, but also involves a less conservative stepsize update rule.
		More generally, we further waive the computation of the norm of the linear operator (required by the stepsize update rule) through a linesearch procedure.
		The resulting algorithm, \refadaPD+*, can be implemented efficiently given that the backtracks do not require additional gradient evaluations.

		Future research directions involve extending the presented line of proof to the setting of variational inequalities in the framework of algorithms such as the extra-gradient method and other projection-based splittings such as those in \cite{latafat2017asymmetric,giselsson2021nonlinear}.
		Other potential directions include nonmonotone extensions for problems satisfying the so-called weak Minty assumptions \cite{diakonikolas2021efficient,pethick2022escaping}, as well as block-coordinate and stochastic variants.
		It would also be interesting to explore extensions of the ideas presented in this work to other types of adaptive settings such as those in \cite{li2019convergence,defazio2022grad,yurtsever2021three,ward2019adagrad}.

	\FloatBarrier
	\begin{appendix}
		\phantomsection\addcontentsline{toc}{section}{Appendix}%

		\section{Convergence analysis for \texorpdfstring{\cref{sec:PD}}{the primal-dual algorithms}}\label{sec:PD:convergence}
			We begin by establishing an inequality for the iterates in \eqref{eq:PD:update} that would eventually under proper stepsize choices lead to a descent inequality.

			\begin{lemma}\label{thm:PD:descent}%
				Let \(\seq{\gamk}, \seq{\sigk}\) be sequences of strictly positive scalars, and starting from a triplet \((x^{-1},x^0,y^0)\in\R^n\times\R^n\times\R^m\) let \(\seq{x^k,y^k}\) be recursively defined by \eqref{eq:PD:update}.
				Then, for any \(\varepsilon_k, \tau_k, \mu_k>0\) and \(\eta_k\geq\frac{\|\linop*(y^k-y^{k-1})\|}{\|y^k-y^{k-1}\|}\) (\eg \(\eta_k=\|\linop\|\)) it holds that
				\begin{align*}
					0
				\leq{} &
					\tfrac{1}{2}\|x^k-x^\star\|^2
					-
					\tfrac{1}{2}\|x^{k+1}-x^\star\|^2
					+
					\tfrac{\varepsilon_{k+1}\rhok*+\mu_{k+1}\eta_{k+1}\gamk*-1}{2}
					\|x^k-x^{k+1}\|^2
				\\
				&
					+
					\rhok*
					\left(
						\tfrac{1-\alpha_k(2-\beta_k)}{2\varepsilon_{k+1}}
						+
						\tfrac{\tau_{k+1}\eta_{k+1}\gamk*}{2}
						-
						\rhok*(1-\alpha_k)
					\right)
					\|x^{k-1}-x^k\|^2
				\\
				&
					+
					\tfrac{\gamk*}{2\sigk*}\|y^k-y^\star\|^2
					-
					\left(
						\tfrac{\gamk*}{2\sigk*}
						-
						\tfrac{\eta_{k+1}\gamk*}{2}\left(
							\tfrac{1}{\mu_{k+1}}
							+
							\tfrac{\rhok*}{\tau_{k+1}}
						\right)
					\right)
					\|y^{k+1}-y^k\|^2
				\\
				&
					-
					\tfrac{\gamk*}{2\sigk*}\|y^{k+1}-y^\star\|^2
				\\
				&
					+
					\gamk*\rhok*P_{k-1}
					-
					\gamk*(1+\rhok*)P_k
					-
					\gamk*Q_{k+1},
				\end{align*}
				where \(\rhok*\coloneqq\frac{\gamk*}{\gamk}\), \(\alpha_k\coloneqq\gamk\Lk\), \(\beta_k\coloneqq\gamk\ck\), and \(P_k, Q_k\) are as in \eqref{subeq:PQ}.
				\begin{proof}
					Let
					\(
						\bar y^{k+1}
					\coloneqq
						y^k+\sigk*\tfrac{\gamk*}{\gamk}\linop (x^k-x^{k-1})
					\),
					so that \(y^{k+1}=\prox_{\sigk*\conj h}(\bar y^{k+1}+\sigk*\linop x^k)\).
					The characterization of \(y^k\) and \(x^{k+1}\) as in the respective updates then reads
					\begin{subequations}\label{eq:subgrad:xy}
						\begin{gather}
						\label{eq:subgrad:xyh}
							\tfrac{\bar y^{k+1}-y^{k+1}}{\sigk*}+\linop x^k
						\in
							\partial \conj h(y^{k+1})
						\shortintertext{and}
						\label{eq:subgrad:xyg}
							\tfrac{\Hk*(x^k)-x^{k+1}}{\gamk*} -\linop*y^{k+1}
						=
							\tfrac{x^k-x^{k+1}}{\gamk*}-\nabla f(x^k)-\linop*y^{k+1}
							\in
							\partial g(x^{k+1}).
						\end{gather}
					\end{subequations}
					In particular, for any solution pair \((x^\star,y^\star)\in \mathcal S_\star\)
					\begin{align*}
						0
					\overrel[\leq]{\eqref{eq:subgrad:xyg}}{} &
						g(x^\star)-g(x^{k+1})
						+
						\innprod{\linop*y^{k+1}}{x^\star-x^{k+1}}
						+
						\overbracket*[0.5pt]{
							\innprod{\nabla f(x^k)}{x^\star-x^{k+1}}
						}^{(\text{A}')}
					\\
					&
						+
						\tfrac{1}{2\gamk*}\|x^k-x^\star\|^2
						-
						\tfrac{1}{2\gamk*}\|x^{k+1}-x^\star\|^2
						-
						\tfrac{1}{2\gamk*}\|x^k-x^{k+1}\|^2
					\shortintertext{and}
						0
					\overrel[\leq]{\eqref{eq:subgrad:xyh}}{} &
						\conj h(y^\star)-\conj h(y^{k+1})
						-
						\innprod{\linop x^k}{y^\star-y^{k+1}}
						-
						\tfrac{1}{\sigk*}
						\innprod{\bar y^{k+1}-y^{k+1}}{y^\star-y^{k+1}}
					\\
					={} &
						\conj h(y^\star)-\conj h(y^{k+1})
						-
						\innprod{\linop u^{k+1}}{y^\star-y^{k+1}}
					\\
					&
						+
						\tfrac{1}{2\sigk*}\|y^k-y^\star\|^2
						-
						\tfrac{1}{2\sigk*}\|y^{k+1}-y^k\|^2
						-
						\tfrac{1}{2\sigk*}\|y^{k+1}-y^\star\|^2,
					\end{align*}
					where \(u^{k+1}\coloneqq (1+\rhok*)x^k-\rhok*x^{k-1}\).
					We next proceed to upper bound the term \((\text{A}')\) as
					\begin{align*}
						(\text{A}')
					={} &
						\innprod{\nabla f(x^k)}{x^\star-x^k}
						+
						\innprod{\nabla f(x^k)}{x^k-x^{k+1}}
					\\
					={} &
						\innprod{\nabla f(x^k)}{x^\star-x^k}
						+
						\tfrac{1}{\gamk}\innprod{\Hk(x^{k-1})-x^k}{x^{k+1}-x^k}
					\\
					&
						+
						\tfrac{1}{\gamk}\innprod{\Hk(x^{k-1})-\Hk(x^k)}{x^k-x^{k+1}}
					\\
					\overrel[\leq]{\eqref{eq:subgrad:xyg}}{} &
						\fillwidthof[c]{\innprod{\nabla f(x^k)}{x^\star-x^k}}{
							f(x^\star) - f(x^k)
						}
						+
						g(x^{k+1})-g(x^k)
						+
						\innprod{\linop*y^k}{x^{k+1}-x^k}
					\\
					&
						+
						\underbracket[0.5pt]{
							\tfrac{1}{\gamk}\innprod{\Hk(x^{k-1})-\Hk(x^k)}{x^k-x^{k+1}}
						}_{\text{(B)}}.
					\end{align*}
					Let \(\varphi = f +g\).
					We bound the term (B) as done in \eqref{eq:B:CL} and combine the three inequalities to obtain
					\begin{align*}
						0
					\leq{} &
						(f+g)(x^\star)-(f+g)(x^k)
						+
						\tfrac{1}{2\gamk*}\|x^k-x^\star\|^2
						-
						\tfrac{1}{2\gamk*}\|x^{k+1}-x^\star\|^2
					\\
					&
						+
						\left(
							\tfrac{\varepsilon_{k+1}}{2\gamk}
							-
							\tfrac{1}{2\gamk*}
						\right)
						\|x^k-x^{k+1}\|^2
						+
						\tfrac{1-\gamk\Lk(2-\gamk\ck)}{2\varepsilon_{k+1}\gamk}\|x^{k-1}-x^k\|^2
					\\
					&
						+
						\conj h(y^\star)-\conj h(y^{k+1})
						+
						\innprod{\linop*y^k}{x^{k+1}-x^k}
						+
						\innprod{\linop*y^{k+1}}{x^\star-x^{k+1}}
						-
						\innprod{\linop u^{k+1}}{y^\star-y^{k+1}}
					\\
					&
						+
						\tfrac{1}{2\sigk*}\|y^k-y^\star\|^2
						-
						\tfrac{1}{2\sigk*}\|y^{k+1}-y^k\|^2
						-
						\tfrac{1}{2\sigk*}\|y^{k+1}-y^\star\|^2.
					\end{align*}
					Using again the subgradient \eqref{eq:subgrad:xy}, one has
					\begin{equation}\label{eq:vk:PD}
						v^k
					\coloneqq
						\tfrac{x^{k-1}-x^k}{\gamk}-(\nabla f(x^{k-1})-\nabla f(x^k))-\linop*y^k
					\in
						\partial(f+g)(x^k),
					\end{equation}
					hence
					\begin{align*}
						0
					\leq{} &
						\tfrac{\gamk*}{\gamk}\left(
							(f+g)(x^{k-1})-(f+g)(x^k)
							-
							\innprod{v^k}{x^{k-1}-x^k}
						\right)
					\\
					={} &
						\tfrac{\gamk*}{\gamk}\left(
							(f+g)(x^{k-1})-(f+g)(x^k)
						\right)
						-
						\tfrac{\gamk*}{\gamk}\tfrac{1-\gamk\Lk}{\gamk}
						\|x^k-x^{k-1}\|^2
					\\
					&
						+
						\tfrac{\gamk*}{\gamk}
						\innprod{\linop*y^k}{x^{k-1}-x^k}.
					\numberthis\label{eq:PD:subgradineq}
					\end{align*}
					Sum the last two inequalities and use the identity \(\rhok*(x^k-x^{k-1})=u^{k+1}-x^k\) to obtain
					\begin{align*}
						0
					\leq{} &
						\tfrac{1}{2\gamk*}\|x^k-x^\star\|^2
						-
						\tfrac{1}{2\gamk*}\|x^{k+1}-x^\star\|^2
						+
						\innprod{\linop x^\star}{y^{k+1}}
						-
						\innprod{\linop u^{k+1}}{y^\star}
					\\
					&
						+
						\left(
							\tfrac{\varepsilon_{k+1}}{2\gamk}
							-
							\tfrac{1}{2\gamk*}
						\right)
						\|x^k-x^{k+1}\|^2
						+
						\left(
							\tfrac{1-\gamk\Lk(2-\gamk\ck)}{2\varepsilon_{k+1}\gamk}
							-
							\tfrac{\gamk*}{\gamk}\tfrac{1-\gamk\Lk}{\gamk}
						\right)
						\|x^{k-1}-x^k\|^2
					\\
					&
						+
						\tfrac{1}{2\sigk*}\|y^k-y^\star\|^2
						-
						\tfrac{1}{2\sigk*}\|y^{k+1}-y^k\|^2
						-
						\tfrac{1}{2\sigk*}\|y^{k+1}-y^\star\|^2
					\\
					&
						+
						\left(\conj h(y^\star)-\conj h(y^{k+1})\right)
					\\
					&
						+
						\underbracket[0.5pt]{%
							\tfrac{\gamk*}{\gamk}\left(
								(f+g)(x^{k-1})-(f+g)(x^k)
							\right)
							-
							\left((f+g)(x^k)-(f+g)(x^\star)\right)
						}_{\text{(C)}}
					\\
					&
						{}+
						\underbracket[0.5pt]{%
							\vphantom{\tfrac{\gamk*}{\gamk}}
							\innprod{\linop*(y^{k+1}-y^k)}{u^{k+1}-x^{k+1}}
						}_{\text{(D)}}.
						\vphantom{\underbracket[0.5pt]{A}_{}}
					\end{align*}
					To conclude, observe that
					\begin{align*}
						\text{(C)}
					={} &
						\rhok*P_{k-1}
						-
						(1+\rhok*)P_k
						+
						\innprod{u^{k+1}-x^\star}{\linop*y^\star}
					\shortintertext{and}
						\text{(D)}
					={} &
						\tfrac{\gamk*}{\gamk}
						\innprod{x^k-x^{k-1}}{\linop*(y^{k+1}-y^k)}
						+
						\innprod{x^k-x^{k+1}}{\linop*(y^{k+1}-y^k)}
					\\
					\leq{} &
						\tfrac{\gamk*}{\gamk}
						\tfrac{\tau_{k+1}\eta_{k+1}}{2}
						\|x^{k-1}-x^k\|^2
						+
						\tfrac{\mu_{k+1}\eta_{k+1}}{2}\|x^{k+1}-x^k\|^2
					\\
					&
						{}+
						\tfrac{\eta_{k+1}}{2}\left(
							\tfrac{1}{\mu_{k+1}}
							+
							\tfrac{\gamk*}{\gamk\tau_{k+1}}
						\right)
						\|y^{k+1}-y^k\|^2,
					\numberthis\label{eq:uLy}
					\end{align*}
					where \(\mu_k, \tau_k>0\) are parameters related to Young's inequality, so that the proof follows from the identity
					\(
						\innprod{\linop x^\star}{y^{k+1}-y^\star}
					=
						\SD(x^\star,y^{k+1})-\SD(x^\star,y^\star)
					\).
				\end{proof}
			\end{lemma}

			In addition to the upper estimate \(L_{f, \V}\) for \(L_k\) and  \(\Lk\) in \cref{thm:Lk<=Ck}, the proof of \cref{thm:PD} will exploit an upper bound for \(\ck\) which is obtained in the next lemma. The proof is a slight refinement of known cocoercivity results in the globally Lipschitz setting \cite[Thm. 5.8]{beck2017first}, and is provided for completeness to highlight the need of the enlarged set \(\overline\V\)
			\begin{lemma}\label{lem:enlarged:ck}
				Suppose that \cref{ass:f} holds.
				Then, for every \(x,y\in\R^n\) it holds that
				\[
					f(y)
				\leq
					f(x)
					+
					\innprod{\nabla f(x)}{y-x}
					+
					\tfrac{1}{2L_{f,\overline\V}}
					\|\nabla f(x) - \nabla f(y)\|^{2},
				\]
				where \(L_{f,\overline{\V}}\) is a Lipschitz modulus for \(\nabla f\) on \(\overline{\V}\coloneqq\V+\cball{0}{\diam(\V)}\), and \(\V\subseteq\R^n\) is a bounded and convex set that contains \(x\) and \(y\).
				In particular,
				\[
					\frac{
						\|\nabla f(x) - \nabla f(y)\|^2
					}{
						\innprod{\nabla f(x) - \nabla f(y)}{x-y}
					}
				\leq
					L_{f,\overline\V}.
				\]
				\begin{proof}
					We henceforth fix \(x,y\in \V\subseteq\overline\V\).
					Since $\nabla f$ is Lipschitz continuous on $\overline\V$ with modulus \(L_{f,\overline\V}\), it follows from the descent inequality \cite[Prop. A.24]{bertsekas2016nonlinear} that
					\begin{equation}\label{eq:holder_ineq}
						f(z)
					\leq
						f(y)
						+
						\innprod{\nabla f(y)}{z-y}
						+
						\tfrac{L_{f,\overline\V}}{2}\|z-y\|^{2}
					\quad
						\forall z \in \overline\V.
					\end{equation}
					Let \(l_x(y) \coloneqq f(y) - f(x) - \innprod{\nabla f(x)}{y-x}\), and note that \(l_x\) is a convex function with \(\nabla l_x(y)=\nabla f(y)-\nabla f(x)\).
					For any \(z \in \overline\V\), we have
					\begin{align*}
						l_x(z)
					={} &
						f(z) - f(x) - \innprod{\nabla f(x)}{z-x}
					\\
					\overrel[\leq]{\eqref{eq:holder_ineq}}{} &
						f(y)
						+
						\innprod{\nabla f(y)}{z-y}
						+
						\tfrac{L_{f,\overline\V}}{2}\|z-y\|^{2}
						-
						f(x)
						-
						\innprod{\nabla f(x)}{z-x}
					\\
					={} &
						f(y) - f(x) - \innprod{\nabla f(x)}{y-x}
						+
						\innprod{\nabla f(y) - \nabla f(x)}{z-y}
						+
						\tfrac{L_{f,\overline\V}}{2}\|z-y\|^{2}
					\\
					={} &
						l_x(y)
						+
						\innprod{\nabla l_x(y)}{z-y}
						+
						\tfrac{L_{f,\overline\V}}{2}\|z-y\|^{2}.
					\end{align*}
					Noticing that $\nabla l_x(x) = 0$, it follows from convexity of \(l_x\) that $x$ is its global minimizer. Hence that \(\min l_x=l_x(x)=0\).
					Let us denote \(v\coloneqq\frac{1}{\|\nabla l_x(y)\|}\nabla l_x(y)\) and set
					\(
						z = y - \frac{\|\nabla l_x(y)\|}{L_{f,\overline\V}}v
					\).
					Note that
					\[
						\|z-y\|
					=
						\tfrac{\|\nabla l_x(y)\|}{L_{f,\overline\V}}
					=
						\tfrac{\|\nabla f(y)-\nabla f(x)\|}{L_{f,\overline\V}}
					\leq
						\|y-x\|,
					\]
					and in particular \(z\in \V+\cball{y}{\|y-x\|}\subseteq\overline\V\).
					From the previous inequality we get
					\begin{qedalign*}
						0
					=
						\min l_x
					=
						l_x(x)
					\leq{} &
						l_x\Bigl(
							y - \tfrac{\|\nabla l_x(y)\|}{L_{f,\overline\V}}v
						\Bigr)
					\\
					\overrel[\leq]{\eqref{eq:holder_ineq}}{} &
						l_x(y)
						-
						\tfrac{\|\nabla l_x(y)\|}{L_{f,\overline\V}}
						\innprod{\nabla l_x(y)}{v}
						+
						\tfrac{L_{f,\overline\V}}{2}
						\tfrac{1}{L_{f,\overline\V}^{2}}
						\|\nabla l_x(y)\|^{2}
					\\
					={} &
						f(y) - f(x) - \innprod{\nabla f(x)}{y-x}
						-
						\tfrac{1}{2L_{f,\overline\V}}
						\|\nabla f(x) - \nabla f(y)\|^{2}.
					\end{qedalign*}
				\end{proof}
			\end{lemma}

			\begin{appendixproof}{thm:PD}[\empty]%
				\begin{proofitemize}
				\item \itemref[e]{thm:PD:SD}~
					Denoting
					\[
						\tilde\U_k
					\coloneqq
						\tfrac{1}{2}\|x^k-x^\star\|^2
						+
						\tfrac{1-\varepsilon_k\rhok-\mu_k\eta_k\gamk}{2}
						\|x^k-x^{k-1}\|^2
						+
						\tfrac{\gamk}{2\sigk}\|y^k-y^\star\|^2
						+
						\gamk(1+\rhok)P_{k-1},
					\]
					the inequality in \cref{thm:PD:descent} can be expressed as
					\begin{align*}
						\tilde\U_{k+1}
					\leq{} &
						\tilde\U_k
						-
						\left(
							\tfrac{\gamk}{2\sigk}
							-
							\tfrac{\gamk*}{2\sigk*}
						\right)
						\|y^k-y^\star\|^2
						-
						\gamk(1+\rhok-\rhok*^2)P_{k-1}
					\\
					&
						-
						\left(
							\tfrac{1-\varepsilon_k\rhok-\mu_k\eta_k\gamk}{2}
							+
							\rhok*
							\left(
								\rhok*(1-\alpha_k)
								-
								\tfrac{1-\alpha_k(2-\beta_k)}{2\varepsilon_{k+1}}
								-
								\tfrac{\tau_{k+1}\eta_{k+1}\gamk*}{2}
							\right)
						\right)
						\|x^{k-1}-x^k\|^2
					\\
					&
						-
						\left(
							\tfrac{\gamk*}{2\sigk*}
							-
							\tfrac{\eta_{k+1}\gamk*}{2}\left(
								\tfrac{1}{\mu_{k+1}}
								+
								\tfrac{\rhok*}{\tau_{k+1}}
							\right)
						\right)
						\|y^{k+1}-y^k\|^2
						-
						\gamk*\bigl(
							\SD(x^\star,y^\star)
							-
							\SD(x^\star,y^{k+1})
						\bigr).
					\end{align*}
					Since \(\sigk=t^2\gamk\) and \(\xik=t^2\gamk^2\eta_k^2\), by selecting \(\varepsilon_k\coloneqq\frac{1}{2\rhok}\) and \(\mu_k\coloneqq2 (1+\sdcoeff)t \xik^{\nicefrac12}\) one has that \(\tilde\U_k=\Uk\).
					Note also that the coefficient of \(\|x^k-x^{k-1}\|\) in \(\Uk\) is strictly positive since the stepsize update ensures \(\xik \leq \nicefrac{1}{4\const^2}\).
					Moreover, with this choice the above inequality becomes
					\begin{align*}
						\Uk*
					\leq{} &
						\Uk
						-
						\gamk(1+\rhok-\rhok*^2)P_{k-1}
						-
						\gamk*\bigl(
							\SD(x^\star,y^\star)
							-
							\SD(x^\star,y^{k+1})
						\bigr)
					\\
					&
						-
						\left(
							\tfrac{1-4\xik(1+\sdcoeff)}{4}
							+
							\rhok*\left(
								\rhok*\alpha_k(1-\beta_k)
								-
								\tfrac{\tau_{k+1}\xik*^{\nicefrac12}}{2t}
							\right)
						\right)
						\|x^{k-1}-x^k\|^2
					\\
					&
						-
						\tfrac{1}{2t}
						\left(
							\tfrac{1+2\sdcoeff}{2t(1+\sdcoeff)}
							-
							\tfrac{\xik*^{\nicefrac12}\rhok*}{\tau_{k+1}}
						\right)
						\|y^{k+1}-y^k\|^2.
					\end{align*}
					We now set \(\tau_{k+1}\coloneqq\rhok*\mu_{k+1}=2(1+\sdcoeff)t\xik*^{\nicefrac12}\rhok*\) so the inequality overall simplifies to the one of the statement
					\begin{align*}
						\Uk*
					\leq{} &
						\Uk
						-
						\gamk(1+\rhok-\rhok*^2)P_{k-1}
						-
						\gamk*\bigl(
							\SD(x^\star,y^\star)
							-
							\SD(x^\star,y^{k+1})
						\bigr)
					\\
					&
						-
						\underbracket[0.5pt]{%
							\left(
								\tfrac{1-4\xik(1+\sdcoeff)}{4}
								+
								\rhok*\left(
									\rhok*\alpha_k(1-\beta_k)
									-
									\rhok*\xik*(1+\sdcoeff)
								\right)
							\right)
						}_{\geq \tfrac{\sdcoeff}{4(1+\sdcoeff)}}
						\|x^{k-1}-x^k\|^2
					\\
					&
						-
						\tfrac{\sdcoeff}{2t^2(1+\sdcoeff)}
						\|y^{k+1}-y^k\|^2,
					\end{align*}
					up to ensuring that the coefficient of \(P_{k-1}\) is positive and that the inequality for the coefficient of \(\|x^k-x^{k-1}\|^2\) holds.
					The former is of trivial verification, having \(\rhok*\leq\sqrt{1+\rhok}\).
					It thus remains to show that
					\[
						\tfrac14
						-
						\tfrac{\sdcoeff}{4(1+\sdcoeff)}
						-
						\xik(1+\sdcoeff)
						-
						\rhok*^2(\delk\mathbin{+}\xik*(1+\sdcoeff))
					\geq
						0
					\]
					holds for every \(k\).
					By using the fact that \(\xik*=(t\eta_{k+1}\gamk*)^2=(t\eta_{k+1}\gamk)^2\rhok*^2\), this reduces to the second-order inequality (in \(\rhok*^2\))
					\[
						(t\eta_{k+1}\gamk)^2(1+\sdcoeff)\rhok*^4
						+
						\delk\rhok*^2
						-
						\bigl[
							\tfrac{1}{4(1+\sdcoeff)}
							-
							\xik(1+\sdcoeff)
						\bigr]
					\leq
						0.
					\]
					Note that the bound
					\(
						\gamk\leq\frac{1}{2\const t\eta_k}
					\)
					implies \(\tfrac1{4(1+\sdcoeff)}-\xik(1+\sdcoeff)\geq\frac{\const^2-(1+\sdcoeff)^2}{4\const^2(1+\sdcoeff)}>0\), thus ensuring that the inequality always admits solutions for small enough \(\rhok*^2\).
					Namely, letting \(\bar\xi_k \coloneqq \xik (1+\sdcoeff)^2\)
					\[
						\rhok*^2
					\leq
						\frac{
							-\delk
							+
							\sqrt{\delk^2+(t\eta_{k+1}\gamk)^2(1-4\bar\xi_k)}
						}{
							2(1+\sdcoeff)(t\eta_{k+1}\gamk)^2
						}
					=
						\frac{
							1-4\bar\xi_k
						}{
							2(1+\sdcoeff)
							\left(
								\delk
								+
								\sqrt{\delk^2+(t\eta_{k+1}\gamk)^2(1-4\bar\xi_k)}
							\right)
						},
					\]
					which is indeed guaranteed by one of the bounds on \(\gamk*=\gamk\rhok*\).
					Note that the second expression removes the singularity in case \(\eta_{k+1}=0\).

				\item \itemref[e]{thm:PD:gamk}~
					Boundedness of the sequence follows from the fact that
					\(
						\tfrac{1}{2}\|x^k-x^\star\|^2
						+
						\tfrac{1}{2t^2}\|y^k-y^\star\|^2
					\leq
						\Uk
					\leq
						\U_1
					\),
					where the first inequality follows by definition of \(\U_k\), cf. \eqref{eq:PD:U}, and the second one from assertion \itemref{thm:PD:SD}.
					In particular, there exists a convex and compact set \(\V\subseteq\R^n\) that contains \(\seq{x^k}\).
					Up to considering a suitable enlargement \(\overline\V\) as in \cref{lem:enlarged:ck},  					\(\Lk\leq\ck\leq L_{f,\overline{\V}}\) holds for every \(k\) by \cref{thm:Lk<=Ck} and \cref{lem:enlarged:ck}, where \(L_{f,\overline\V}>0\) is a Lipschitz modulus for \(\nabla f\) on \(\overline\V\).
					To prove the lower bound on the stepsize, we will show that whenever \(\gamk*<\gamk\) occurs, then necessarily \(\gamk*\) is greater than some constant \(\hat\gamma\) as in the statement.
					The proof will then follow from a trivial inductive argument.
					Suppose that \(\gamk*<\gamk\).
					If \(\gamk*=\frac{1}{2\const t\eta_{k+1}}\), then clearly \(\gamk*\geq\frac{1}{2\const t\eta_{\rm max}}\).
					Otherwise, necessarily
					\begin{align}
					\nonumber
						\gamk^2
					>
						\gamk*^2
					={} &
						\frac{
							\gamk^2(1-4\bar\xi_k)
						}{
							2(1+\sdcoeff)
							\left(
								\sqrt{
									\delk^2
									+
									(t\eta_{k+1}\gamk)^2(1-4\bar\xi_k)
								}
								+
								\delk
							\right)
						}
					\\
					\geq{} &
						\frac{
							\gamk^2 \bar\const
						}{
							2(1+\sdcoeff)
							\left(
								\sqrt{
									\delk^2
									+
									(t\eta_{\rm max}\gamk)^2 \bar\const
								}
								+
								\delk
							\right)
						},
						\numberthis\label{eq:PD:gamkdecreaseDelta}
					\end{align}
					where
					\(
						\bar\const
					\coloneqq
						\frac{\const^2-(1 + \sdcoeff)^2}{\const^2}
					\),
					and the second inequality uses the fact that \(1-4\bar\xi_k\geq \bar\const>0\) together with the fact that \((0,\infty)\ni x\mapsto\frac{x}{\sqrt{b^2+a^2x\,}+b}\) is increasing for any value of \(a,b\in\R\).
					We now distinguish two cases:
					\begin{proofitemize}
					\item Case 1.
						\(\delk\leq0\).
						By comparing the outermost terms of the chain of inequalities in \eqref{eq:PD:gamkdecreaseDelta} we obtain that
						\[
							\sqrt{
								\delk^2
								+
								(t\eta_{\rm max}\gamk)^2\bar\const
							}
							+
							\delk
						>
							\tfrac{\bar\const}{2(1+\sdcoeff)}
						\quad\Leftrightarrow\quad
							\delk
						>
							\tfrac{\bar\const}{4(1+\sdcoeff)}-(1+\sdcoeff)(t\eta_{\rm max}\gamk)^2,
						\]
						and in particular \(\gamk>\frac{\sqrt{\bar\const}}{2(1+\sdcoeff)t\eta_{\rm max}}\).
						Then, since \(x\mapsto\frac{1}{\sqrt{x^2+b^2\,}+x}\) is decreasing for any value of \(b\in\R\), by setting \(\delk\) equal to 0 in \eqref{eq:PD:gamkdecreaseDelta} we obtain
						\[
							\gamk*^2
						\geq
							\gamk
							\frac{
								\sqrt{\bar\const}
							}{
								2(1+\sdcoeff)t\eta_{\rm max}
							}
						\geq
							\frac{
								\bar\const
							}{
								(2(1+\sdcoeff)t\eta_{\rm max})^2
							}.
						\]

					\item Case 2.
						\(\delk>0\) or, equivalently, \(\gamk\ck>1\).
						Denoting \(\alpha\coloneqq\gamk L_{f,\overline{\V}}\), one has that \(\alpha\geq\gamk\Lk\) and \(\alpha\geq\gamk\ck>1\), hence that \(\delk\leq\alpha(\alpha-1)\).
						Arguing as in the previous case, this time by setting \(\delk\gets\alpha(\alpha-1)\) in \eqref{eq:PD:gamkdecreaseDelta}, yields
						\begin{align*}
							\gamk*^2L_{f,\overline{\V}}^2
						\geq{} &
							\frac{
								\alpha^2 \bar\const%
							}{
								2(1+\sdcoeff)
								\left(
									\sqrt{
										\alpha^2(\alpha-1)^2
										+
										(t\eta_{\rm max}\gamk)^2 \bar\const %
									}
									+
									\alpha(\alpha-1)
								\right)
							}
						\\
						={} &
							\frac{
								\alpha \bar\const
							}{
								2(1+\sdcoeff)
								\left(
									\sqrt{
										(\alpha-1)^2
										+
										(\nicefrac{t\eta_{\rm max}}{L_{f,\overline{\V}}})^2 \bar\const  %
									}
									+
									\alpha-1
								\right)
							}
						\\
						\geq{} &
							\min\set{
								\tfrac{
									L_{f,\overline{\V}}\sqrt{\bar\const}
								}{
									2(1+\sdcoeff)
									t\eta_{\rm max}
								},\,
								\tfrac{
									\bar\const
								}{
									4(1+\sdcoeff)
								}
							},
						\end{align*}
						where the last inequality owes to the fact that \([1,\infty)\ni x\mapsto\frac{x}{\sqrt{(x-1)^2+b^2\,}+x-1}\) attains the infimum at either \(1\) or \(\infty\) for any \(b\in\R\).
					\end{proofitemize}
					Putting all the cases together yields
					\begin{equation}\label{eq:PDgamkLB}
						\gamk
					\geq
						\hat\gamma
					\coloneqq
						\min\set{
							\gamma_0
						,\,
							\tfrac{
								\sqrt[4]{\bar\const}
							}{
								\sqrt{
								2(1+\sdcoeff)
								t\eta_{\rm max}L_{f,\overline{\V}}
								}
							}
						,\,
							\tfrac{
								\sqrt{\bar\const}
							}{
								2\sqrt{(1+\sdcoeff)}L_{f,\overline{\V}}
							}
						,\,
							\tfrac{1}{2\const t\eta_{\rm max}}
						,\,
							\tfrac{
									\sqrt{\bar\const}
								}{
									2(1+\sdcoeff)t\eta_{\rm max}
							}
						}
					>
						0
					\end{equation}
					(where we remind that \(\bar\const=\frac{\const^2-(1 + \sdcoeff)^2}{\const^2}\)), establishing the claim.

				\item \itemref[e]{thm:PD:PQ}~
					That \(Q_k\to0\) as \(k\to\infty\) follows from a telescoping argument in assertion \itemref[e]{thm:PD:SD}, since \(\seq{\gamk}\) is bounded away from zero by assertion \itemref[e]{thm:PD:gamk}.
					As to \(\seq{P_k}\), if \(\limsup_{k\to\infty}(1+\rhok-\rhok*^2)>0\), then the same telescoping argument yields the claim.
					Alternatively, by the stepsize update rule one has that \(1+\rhok-\rhok*^2\geq0\) and therefore \(1+\rhok-\rhok*^2\to0\), from which it easily follows that \(\liminf_{k\to\infty}\rhok>1\) and that therefore \(\gamk\to\infty\).
					Consequently, since \(\gamk(1+\rhok)P_{k-1}\leq\Uk\leq\U_0\), this directly proves that \(P_k\to0\).

				\item \itemref[e]{thm:PD:z*}~
					Let \((\hat x, \hat y)\) denote a limit point of \(\seq{x^k, y^k}\).
					Telescoping the inequality in assertion \itemref[e]{thm:PD:SD} yields that both \(\|x^{k+1} - x^k\|\) and \(\|y^{k-1} - y^k\|\) vanish.
					By the optimality condition for \cref{state:PD:x+} we have
					\begin{subequations}\label{subeq:PDoptsubseq}%
						\begin{equation}
							\tfrac{1}{\gamk*}(x^k-x^{k+1})
						\in
							\nabla f(x^k) + \linop*y^{k+1} + \partial g(x^{k+1}).
						\end{equation}
						Passing to the limit along the subsequence converging to \((\hat x, \hat y)\), using the fact that \(\seq{\gamk}\) is bounded away from zero by assertion \itemref[e]{thm:PD:gamk} and outer semicontinuity of \(\partial g\) yield \(0\in \nabla f(\hat x) + \partial g(\hat x) + \linop*\hat y\).
						Similarly, for the dual variable,
						\begin{equation}
							\tfrac{1}{\sigma_{k+1}}(y^k-y^{k+1})
							+
							\tfrac{\gamk*}{\gamk}\linop(x^k-x^{k-1})
						\in
							\partial h^*(y^{k+1})-\linop x^k.
						\end{equation}
					\end{subequations}
					A trivial induction argument reveals that \(\tfrac{\gamk*}{\gamk}\) is upper bounded by \(\max\set{\tfrac{\gamma_0}{\gamma_{-1}}, \tfrac12(1+\sqrt{5})}\).
					Therefore, passing to the limit along the same subsequence and recalling that \(\sigk*=t^2\gamk*\) is bounded away from zero yield \(0\in\partial h^*(\hat y)-\linop\hat x\).
					Along with the previous inclusion, primal-dual optimality of the limit pair \((\hat x, \hat y)\) is established.
					Therefore, any limit point of \(\seq{x^k,y^k}\) is a primal-dual optimal pair.
					Suppose that \((x^\infty, y^\infty)\) and \((\hat x^\infty, \hat y^\infty)\) are two primal-dual optimal limit points of \(\seq{x^k, y^k}\).
					Define \(v^k = (x^k, \tfrac1t y^k)\), and consistently \(v^\infty = (x^\infty, \tfrac1t y^\infty)\), \(\hat v^\infty = (\hat x^\infty, \tfrac1t \hat y^\infty)\).
					Observe that
					\[
						\innprod{v^k}{v^\infty - \hat v^\infty}
					=
						\U_k(\hat x^\infty, \hat y^\infty) - \U_k(x^\infty, y^\infty) + \tfrac12\|v^\infty\|^2 - \tfrac12\|\hat v^\infty\|^2,
					\]
					and since \(\seq{\U_k(x^\star, y^\star)}\) is convergent for all primal-dual optimal pairs \((x^\star, y^\star)\), then so is \(\seq{\innprod{v^k}{v^\infty - \hat v^\infty}}\).
					Passing to the limit along the two converging subsequences thus yields
					\(
						\innprod{v^\infty}{v^\infty - \hat v^\infty}
					=
						\innprod{\hat v^\infty}{v^\infty - \hat v^\infty}
					\),
					which after rearranging results in \(\|v^\infty - \hat v^\infty\|^2=0\) establishing uniqueness of the optimal limit point.
				\qedhere
				\end{proofitemize}
			\end{appendixproof}
	\end{appendix}

	\phantomsection
	\addcontentsline{toc}{section}{References}
	\bibliographystyle{plain}
	\bibliography{Bibliography.bib}

\begin{thebibliography}{10}

\bibitem{alacaoglu2023beyond}
Ahmet Alacaoglu, Axel B{\"o}hm, and Yura Malitsky.
\newblock Beyond the golden ratio for variational inequality algorithms.
\newblock {\em Journal of Machine Learning Research}, 24(172):1--33, 2023.

\bibitem{altschuler2023acceleration}
Jason~M. Altschuler and Pablo~A. Parrilo.
\newblock Acceleration by stepsize hedging {II}: Silver stepsize schedule for
  smooth convex optimization.
\newblock {\em arXiv preprint arXiv:2309.16530}, 2023.

\bibitem{attouch2023fast}
H\'edy Attouch, Radu~Ioan Bo{\c{t}}, and Dang-Khoa Nguyen.
\newblock Fast convex optimization via closed-loop time scaling of gradient
  dynamics.
\newblock {\em arXiv preprint arXiv:2301.00701}, 2023.

\bibitem{baillon1977quelques}
Jean-Bernard Baillon and Georges Haddad.
\newblock Quelques propri{\'e}t{\'e}s des op{\'e}rateurs angle-born{\'e}s et
  $n$-cycliquement monotones.
\newblock {\em Israel Journal of Mathematics}, 26(2):137--150, 1977.

\bibitem{barzilai1988two}
Jonathan Barzilai and Jonathan~M Borwein.
\newblock Two-point step size gradient methods.
\newblock {\em IMA Journal of Numerical Analysis}, 8(1):141--148, 1988.

\bibitem{bauschke2017convex}
Heinz~H. Bauschke and Patrick~L. Combettes.
\newblock {\em Convex Analysis and Monotone Operator Theory in {H}ilbert
  Spaces}.
\newblock CMS Books in Mathematics. Springer, 2017.

\bibitem{beck2017first}
Amir Beck.
\newblock {\em First-Order Methods in Optimization}.
\newblock SIAM, Philadelphia, PA, 2017.

\bibitem{bertsekas2016nonlinear}
Dimitri~P. Bertsekas.
\newblock {\em Nonlinear Programming}.
\newblock Athena Scientific, 2016.

\bibitem{bianchi2014primal}
Pascal Bianchi and Walid Hachem.
\newblock A primal-dual algorithm for distributed optimization.
\newblock In {\em IEEE 53rd Annual Conference on Decision and Control (CDC)},
  pages 4240--4245, dec 2014.

\bibitem{bohm2022solving}
Axel B{\"o}hm.
\newblock Solving nonconvex-nonconcave min-max problems exhibiting weak {M}inty
  solutions.
\newblock {\em arXiv preprint arXiv:2201.12247}, 2022.

\bibitem{bot2013douglas}
Radu~Ioan Bo{\c t} and Christopher Hendrich.
\newblock A {D}ouglas-{R}achford type primal-dual method for solving inclusions
  with mixtures of composite and parallel-sum type monotone operators.
\newblock {\em SIAM Journal on Optimization}, 23(4):2541--2565, 2013.

\bibitem{bot2023relaxed}
Radu~Ioan Bo{\c{t}}, Michael Sedlmayer, and Phan~Tu Vuong.
\newblock A relaxed inertial forward-backward-forward algorithm for solving
  monotone inclusions with application to {GAN}s.
\newblock {\em Journal of Machine Learning Research}, 24:1--37, 2023.

\bibitem{bricenoarias2018forward}
Luis~M. Brice{\~n}o-Arias and Damek Davis.
\newblock Forward-backward-half forward algorithm for solving monotone
  inclusions.
\newblock {\em SIAM Journal on Optimization}, 28(4):2839--2871, 2018.

\bibitem{burdakov2019stabilized}
Oleg Burdakov, Yu-Hong Dai, and Na~Huang.
\newblock Stabilized {B}arzilai-{B}orwein method.
\newblock {\em arXiv preprint arXiv:1907.06409}, 2019.

\bibitem{chambolle2011first}
Antonin Chambolle and Thomas Pock.
\newblock A first-order primal-dual algorithm for convex problems with
  applications to imaging.
\newblock {\em Journal of Mathematical Imaging and Vision}, 40(1):120--145,
  2011.

\bibitem{chang2011libsvm}
Chih-Chung Chang and Chih-Jen Lin.
\newblock {LIBSVM}: A library for support vector machines.
\newblock {\em ACM Transactions on Intelligent Systems and Technology (TIST)},
  2:1--27, 2011.

\bibitem{chang2022golden}
Xiao-Kai Chang, Junfeng Yang, and Hongchao Zhang.
\newblock Golden ratio primal-dual algorithm with linesearch.
\newblock {\em SIAM Journal on Optimization}, 32(3):1584--1613, 2022.

\bibitem{combettes2011proximal}
Patrick~L. Combettes and Jean-Christophe Pesquet.
\newblock Proximal splitting methods in signal processing.
\newblock In {\em Fixed-point algorithms for inverse problems in science and
  engineering}, pages 185--212. Springer New York, 2011.

\bibitem{combettes2012primal}
Patrick~L. Combettes and Jean-Christophe Pesquet.
\newblock Primal-dual splitting algorithm for solving inclusions with mixtures
  of composite, {L}ipschitzian, and parallel-sum type monotone operators.
\newblock {\em Set-Valued and variational analysis}, 20(2):307--330, 2012.

\bibitem{condat2013primal}
Laurent Condat.
\newblock A primal-dual splitting method for convex optimization involving
  {L}ipschitzian, proximable and linear composite terms.
\newblock {\em Journal of Optimization Theory and Applications},
  158(2):460--479, 2013.

\bibitem{dai2005projected}
Yu-Hong Dai and Roger Fletcher.
\newblock Projected {B}arzilai-{B}orwein methods for large-scale
  box-constrained quadratic programming.
\newblock {\em Numerische Mathematik}, 100(1):21--47, 2005.

\bibitem{davis2017three}
Damek Davis and Wotao Yin.
\newblock A three-operator splitting scheme and its optimization applications.
\newblock {\em Set-Valued and Variational Analysis}, 25(4):829--858, dec 2017.

\bibitem{demarchi2022proximal}
Alberto De~Marchi and Andreas Themelis.
\newblock Proximal gradient algorithms under local {L}ipschitz gradient
  continuity: A convergence and robustness analysis of {PANOC}.
\newblock {\em Journal of Optimization Theory and Applications}, 194:771--794,
  2022.

\bibitem{defazio2022grad}
Aaron Defazio, Baoyu Zhou, and Lin Xiao.
\newblock Grad-{G}rada{G}rad? {A} non-monotone adaptive stochastic gradient
  method.
\newblock {\em arXiv preprint arXiv:2206.06900}, 2022.

\bibitem{diakonikolas2021efficient}
Jelena Diakonikolas, Constantinos Daskalakis, and Michael Jordan.
\newblock Efficient methods for structured nonconvex-nonconcave min-max
  optimization.
\newblock In {\em International Conference on Artificial Intelligence and
  Statistics}, pages 2746--2754. PMLR, 2021.

\bibitem{drori2015simple}
Yoel Drori, Shoham Sabach, and Marc Teboulle.
\newblock A simple algorithm for a class of nonsmooth convex-concave
  saddle-point problems.
\newblock {\em Operations Research Letters}, 43(2):209--214, 2015.

\bibitem{fercoq2019coordinate}
Olivier Fercoq and Pascal Bianchi.
\newblock A coordinate-descent primal-dual algorithm with large step size and
  possibly nonseparable functions.
\newblock {\em SIAM Journal on Optimization}, 29(1):100--134, 2019.

\bibitem{giselsson2021nonlinear}
Pontus Giselsson.
\newblock Nonlinear forward-backward splitting with projection correction.
\newblock {\em SIAM Journal on Optimization}, 31(3):2199--2226, 2021.

\bibitem{goldstein2015adaptive}
Tom Goldstein, Min Li, and Xiaoming Yuan.
\newblock Adaptive primal-dual splitting methods for statistical learning and
  image processing.
\newblock {\em Advances in neural information processing systems}, 28, 2015.

\bibitem{goldstein2013adaptive}
Tom Goldstein, Min Li, Xiaoming Yuan, Ernie Esser, and Richard Baraniuk.
\newblock Adaptive primal-dual hybrid gradient methods for saddle-point
  problems.
\newblock {\em arXiv preprint arXiv:1305.0546}, 2013.

\bibitem{grimmer2023accelerated}
Benjamin Grimmer, Kevin Shu, and Alex~L Wang.
\newblock Accelerated gradient descent via long steps.
\newblock {\em arXiv preprint arXiv:2309.09961}, 2023.

\bibitem{hastie2001elements}
Trevor Hastie, Jerome Friedman, and Robert Tibshirani.
\newblock {\em The Elements of Statistical Learning}.
\newblock Springer New York, 2001.

\bibitem{he2012convergence}
Bingsheng He and Xiaoming Yuan.
\newblock Convergence analysis of primal-dual algorithms for a saddle-point
  problem: from contraction perspective.
\newblock {\em SIAM Journal on Imaging Sciences}, 5(1):119--149, 2012.

\bibitem{jezierska2012primal}
Anna Jezierska, Emilie Chouzenoux, Jean-Christophe Pesquet, and Hugues Talbot.
\newblock A primal-dual proximal splitting approach for restoring data
  corrupted with {P}oisson-{G}aussian noise.
\newblock In {\em 2012 IEEE International Conference on Acoustics, Speech and
  Signal Processing (ICASSP)}, pages 1085--1088. IEEE, 2012.

\bibitem{komodakis2015playing}
Nikos Komodakis and Jean-Christophe Pesquet.
\newblock Playing with duality: An overview of recent primal-dual approaches
  for solving large-scale optimization problems.
\newblock {\em IEEE Signal Processing Magazine}, 32(6):31--54, nov 2015.

\bibitem{latafat2020distributed}
Puya Latafat.
\newblock {\em Distributed Proximal Algorithms for Large-Scale Structured
  Optimization}.
\newblock PhD thesis, KU Leuven, jul 2020.

\bibitem{latafat2018plug}
Puya Latafat, Alberto Bemporad, and Panagiotis Patrinos.
\newblock Plug and play distributed model predictive control with dynamic
  coupling: A randomized primal-dual proximal algorithm.
\newblock In {\em European Control Conference (ECC)}, pages 1160--1165, jun
  2018.

\bibitem{latafat2019new}
Puya Latafat, Nikolaos~M. Freris, and Panagiotis Patrinos.
\newblock A new randomized block-coordinate primal-dual proximal algorithm for
  distributed optimization.
\newblock {\em IEEE Transactions on Automatic Control}, 64(10):4050--4065, oct
  2019.

\bibitem{latafat2017asymmetric}
Puya Latafat and Panagiotis Patrinos.
\newblock Asymmetric forward--backward--adjoint splitting for solving monotone
  inclusions involving three operators.
\newblock {\em Computational Optimization and Applications}, 68(1):57--93, sep
  2017.

\bibitem{latafat2018primal}
Puya Latafat and Panagiotis Patrinos.
\newblock Primal-dual proximal algorithms for structured convex optimization: A
  unifying framework.
\newblock In Pontus Giselsson and Anders Rantzer, editors, {\em Large-Scale and
  Distributed Optimization}, pages 97--120. Springer International Publishing,
  2018.

\bibitem{latafat2016new}
Puya Latafat, Lorenzo Stella, and Panagiotis Patrinos.
\newblock New primal-dual proximal algorithm for distributed optimization.
\newblock In {\em 55th IEEE Conference on Decision and Control (CDC)}, pages
  1959--1964, dec 2016.

\bibitem{latafat2022bregman}
Puya Latafat, Andreas Themelis, Masoud Ahookhosh, and Panagiotis Patrinos.
\newblock {B}regman {F}inito/\allowbreak{MISO} for nonconvex regularized finite
  sum minimization without {L}ipschitz gradient continuity.
\newblock {\em SIAM Journal on Optimization}, 32(3):2230--2262, 2022.

\bibitem{latafat2022block}
Puya Latafat, Andreas Themelis, and Panagiotis Patrinos.
\newblock Block-coordinate and incremental aggregated proximal gradient methods
  for nonsmooth nonconvex problems.
\newblock {\em Mathematical Programming}, 193(1):195--224, 2022.

\bibitem{latafat2023adabim}
Puya Latafat, Andreas Themelis, Silvia Villa, and Panagiotis Patrinos.
\newblock On the convergence of proximal gradient methods for convex simple
  bilevel optimization.
\newblock {\em arXiv preprint arXiv:2305.03559}, 2023.

\bibitem{li2019convergence}
Xiaoyu Li and Francesco Orabona.
\newblock On the convergence of stochastic gradient descent with adaptive
  stepsizes.
\newblock In {\em The 22nd International Conference on Artificial Intelligence
  and Statistics}, pages 983--992. PMLR, 2019.

\bibitem{malitsky2020golden}
Yura Malitsky.
\newblock Golden ratio algorithms for variational inequalities.
\newblock {\em Mathematical Programming}, 184(1):383--410, 2020.

\bibitem{malitsky2020adaptive}
Yura Malitsky and Konstantin Mishchenko.
\newblock Adaptive gradient descent without descent.
\newblock In {\em Proceedings of the 37th International Conference on Machine
  Learning}, volume 119, pages 6702--6712. PMLR, 13- 2020.

\bibitem{malitsky2023adaptive}
Yura Malitsky and Konstantin Mishchenko.
\newblock Adaptive proximal gradient method for convex optimization.
\newblock {\em arXiv preprint arXiv:2308.02261}, 2023.

\bibitem{malitsky2018first}
Yura Malitsky and Thomas Pock.
\newblock A first-order primal-dual algorithm with linesearch.
\newblock {\em SIAM Journal on Optimization}, 28(1):411--432, 2018.

\bibitem{malitsky2020forward}
Yura Malitsky and Matthew~K. Tam.
\newblock A forward-backward splitting method for monotone inclusions without
  cocoercivity.
\newblock {\em SIAM Journal on Optimization}, 30(2):1451--1472, 2020.

\bibitem{marumo2022parameter}
Naoki Marumo and Akiko Takeda.
\newblock Parameter-free accelerated gradient descent for nonconvex
  minimization, 2022.

\bibitem{nesterov2013gradient}
Yurii Nesterov.
\newblock Gradient methods for minimizing composite functions.
\newblock {\em Mathematical Programming}, 140(1):125--161, aug 2013.

\bibitem{nesterov2006cubic}
Yurii Nesterov and Boris~T. Polyak.
\newblock Cubic regularization of {N}ewton method and its global performance.
\newblock {\em Mathematical Programming}, 108(1):177--205, 2006.

\bibitem{pedregosa2018adaptive}
Fabian Pedregosa and Gauthier Gidel.
\newblock Adaptive three operator splitting.
\newblock In {\em International Conference on Machine Learning}, pages
  4085--4094. PMLR, 2018.

\bibitem{pethick2022escaping}
Thomas Pethick, Puya Latafat, Panagiotis Patrinos, Olivier Fercoq, and Volkan
  Cevher.
\newblock Escaping limit cycles: Global convergence for constrained
  nonconvex-nonconcave minimax problems.
\newblock In {\em International Conference on Learning Representations}, 2022.

\bibitem{raydan1993barzilai}
Marcos Raydan.
\newblock On the {B}arzilai and {B}orwein choice of steplength for the gradient
  method.
\newblock {\em IMA Journal of Numerical Analysis}, 13(3):321--326, 1993.

\bibitem{rockafellar2009variational}
R.~Tyrrell Rockafellar and Roger J.-B. Wets.
\newblock {\em Variational Analysis}, volume 317.
\newblock Springer, 2009.

\bibitem{rockafellar1970convex}
Ralph~T. Rockafellar.
\newblock {\em Convex Analysis}.
\newblock Princeton University Press, 1970.

\bibitem{ryu2020finding}
Ernest~K. Ryu and B\`{\u a}ng~C. V{\~u}.
\newblock Finding the forward-{D}ouglas--{R}achford-forward method.
\newblock {\em Journal of Optimization Theory and Applications},
  184(3):858--876, mar 2020.

\bibitem{salzo2017variable}
Saverio Salzo.
\newblock The variable metric forward-backward splitting algorithm under mild
  differentiability assumptions.
\newblock {\em SIAM Journal on Optimization}, 27(4):2153--2181, 2017.

\bibitem{Sra2012optimization}
Suvrit Sra, Sebastian Nowozin, and Stephen~J Wright.
\newblock {\em Optimization for Machine Learning}.
\newblock MIT Press, 2012.

\bibitem{tan2016barzilai}
Conghui Tan, Shiqian Ma, Yu-Hong Dai, and Yuqiu Qian.
\newblock {B}arzilai-{B}orwein step size for stochastic gradient descent.
\newblock {\em Advances in neural information processing systems}, 29, 2016.

\bibitem{teboulle2022elementary}
Marc Teboulle and Yakov Vaisbourd.
\newblock An elementary approach to tight worst case complexity analysis of
  gradient based methods.
\newblock {\em Mathematical Programming}, pages 1--34, 2022.

\bibitem{thong2020self}
Duong~Viet Thong, Dang Van~Hieu, and Themistocles~M Rassias.
\newblock Self adaptive inertial subgradient extragradient algorithms for
  solving pseudomonotone variational inequality problems.
\newblock {\em Optimization Letters}, 14(1):115--144, 2020.

\bibitem{vladarean2021first}
Maria-Luiza Vladarean, Yura Malitsky, and Volkan Cevher.
\newblock A first-order primal-dual method with adaptivity to local smoothness.
\newblock {\em Advances in Neural Information Processing Systems},
  34:6171--6182, 2021.

\bibitem{vu2013splitting}
B\`{\u a}ng~C. V{\~u}.
\newblock A splitting algorithm for dual monotone inclusions involving
  cocoercive operators.
\newblock {\em Advances in Computational Mathematics}, 38(3):667--681, 2013.

\bibitem{ward2019adagrad}
Rachel Ward, Xiaoxia Wu, and Leon Bottou.
\newblock {A}da{G}rad stepsizes: Sharp convergence over nonconvex landscapes.
\newblock In Kamalika Chaudhuri and Ruslan Salakhutdinov, editors, {\em
  Proceedings of the 36th International Conference on Machine Learning},
  volume~97 of {\em Proceedings of Machine Learning Research}, pages
  6677--6686. PMLR, 09- 2019.

\bibitem{yan2018new}
Ming Yan.
\newblock A new primal--dual algorithm for minimizing the sum of three
  functions with a linear operator.
\newblock {\em Journal of Scientific Computing}, 76(3):1698--1717, 2018.

\bibitem{Yang2021Self}
Jun Yang.
\newblock Self-adaptive inertial subgradient extragradient algorithm for
  solving pseudomonotone variational inequalities.
\newblock {\em Applicable Analysis}, 100(5):1067--1078, 2021.

\bibitem{yang2018modified}
Jun Yang and Hongwei Liu.
\newblock A modified projected gradient method for monotone variational
  inequalities.
\newblock {\em Journal of Optimization Theory and Applications},
  179(1):197--211, 2018.

\bibitem{yurtsever2021three}
Alp Yurtsever, Alex Gu, and Suvrit Sra.
\newblock Three operator splitting with subgradients, stochastic gradients, and
  adaptive learning rates.
\newblock {\em Advances in Neural Information Processing Systems},
  34:19743--19756, 2021.

\end{thebibliography}

\end{document}